\documentclass[preprint,12pt,3p,sort&compress]{elsarticle}

\usepackage{csquotes}
\usepackage{caption}
\usepackage{ upgreek }
\usepackage{subcaption}
\usepackage{easyReview}
\usepackage{natbib}
\usepackage{amssymb,bm}
\usepackage[mathcal]{euscript}
\usepackage{epstopdf}
\usepackage{amsthm,amsmath,mathtools}

\usepackage{color}
\usepackage{mathrsfs}
\usepackage{dsfont}
\usepackage{graphicx}
\usepackage[shortlabels]{enumitem}
\usepackage{empheq}
\usepackage{float}
\usepackage{physics}
\usepackage[linkcolor=blue, urlcolor=blue, citecolor=blue,colorlinks, bookmarks]{hyperref}
\usepackage{amsfonts}
\usepackage{multirow,booktabs}
\usepackage{adjustbox}
\usepackage{cases} 
\usepackage{algorithm}       
\usepackage{algpseudocode}
\numberwithin{equation}{section}

\numberwithin{Example}{section}
\usepackage{array}
\newcolumntype{L}{>{$}l<{$}} 
\newcolumntype{C}{>{$}c<{$}}
\definecolor{otherblue}{rgb}{0,0.3,0.6}

\begin{document}	
	
	\begin{frontmatter}
		\title{{DAS-PINNs for  high-dimensional partial differential equations: extending deep adaptive sampling to  spacetime domains}}
		\author{Anshima Singh\fnref{label4}}
		 \fntext[label4]{Corresponding author}
		\ead{anshima.singh@manchester.ac.uk}
		\author{David J. Silvester}
		\address{Department of Mathematics, University of Manchester, Oxford Road, Manchester M13 9PL, UK}
		\ead{david.silvester@manchester.ac.uk}
		
\begin{abstract}
{
Time-dependent high-dimensional partial differential equations (PDEs) with spatially localised and dynamically evolving solutions pose a fundamental challenge for physics-informed neural networks (PINNs), as uniform collocation sampling becomes increasingly ineffective in high-dimensional spatiotemporal domains. In this work, a deep adaptive sampling framework for PINNs is extended to the time-dependent setting by treating space and time as a unified domain without any explicit time marching. A normalising flow neural network model effectively learns the distribution induced by the PDE residual and generates new collocation points concentrated in regions where the solution is most difficult to learn. Unlike conventional adaptive strategies that require explicit time stepping or moving meshes, high-residual regions are automatically identified and tracked across both space and time, driven purely by the PDE residual distribution. The effectiveness of the proposed strategy is  assessed on a range of  benchmark problems, from sharp and moving features in two spatial dimensions to localised structures in up to eight spatial dimensions.}
\end{abstract}
		
	\end{frontmatter}
	
\section{Introduction}\label{sec:Intro}
	
{
The approximation of solutions of  partial differential equations (PDEs)  is fundamentally important from a human perspective. Computational simulations of mathematical models of physical phenomena  are essential for applications ranging  from predicting the weather  to  the design of  safe road bridges and  efficient mobile phone networks~\cite{ESW2015}.  Computational fluid dynamics enables the modelling of the flow of blood in the human heart, making the provision  of  patient-specific surgery possible.  Modelling of financial products such as derivatives has  changed the face of  stock market  modelling, but involves the  solution of high-dimensional  PDEs.   Computational chemists also need to  solve high-dimensional PDE problems if they are to model realistic quantum mechanical systems. Solutions to high-dimensional problems are  difficult to approximate  using conventional  discretisation strategies based on finite elements or finite differences. This has opened the door to machine learning strategies such as physics inspired neural networks (PINNS), see~\cite{MR3881695} and neural network approximation of  backward stable stochastic  differential equation (deep BSDEs)~\cite{han2018}, both of which have the  potential to overcome the curse of dimensionality.
}

{
We  explore  adaptive sampling for physics inspired neural  networks in this paper. The basic strategy (DAS--PINNs) was developed for steady-state problems (elliptic PDEs) by Tang et al., see~\cite{MR4531552}. The methodology was subsequently  extended to parametric ODES and steady-state parametric PDEs by Wang et al.~\cite{wang2024}. Our contribution complements the recent work reported by Zhai et al.~\cite{xu2026moving}, wherein localised features are tracked by decoupling the dynamics of the PDEs from the generation of adaptive samples in space. The decoupled strategy in~\cite{xu2026moving} is referred to as a {\sl moving sample method} and is effected  by monitoring the dynamics of adaptive collocation points based on the time varying residual of the neural network solution and then using it for generating samples  at the next adaptive stage. }

{
The main contribution of this work is to show that efficient and reliable approximation  of spatially localised features moving
in time can be realised without decoupling the dynamics from the adaptivity.  (Several of benchmark problems  are taken directly from~\cite{xu2026moving} to facilitate a direct comparison.) Our intuition is that moving mesh or point strategies are  well suited to the solution of  pure transport problems, but are less well suited to parabolic PDE problems  where  the localised features 
change  over time due to diffusive effects.  
}

The remainder of the paper is organised as follows. Section~\ref{sec:Preliminaries} introduces the problem formulation and the PINNs framework. The deep adaptive sampling methodology for spatiotemporal domains is presented in Section~\ref{sec:DAS-PINNs}. Numerical experiments on {a selection of  low-dimensional}  and high-dimensional benchmark problems are reported in Section~\ref{sec:Results}, and the paper is concluded in Section~\ref{sec:Conclusion}.

\section{Preliminaries}\label{sec:Preliminaries}

We  {focus attention on a  general  semilinear,  parabolic or hyperbolic, second-order partial differential equation in this work. The problem specification  is as follows:}
\begin{equation}
	\begin{cases}
	\mathcal{L}u(\mathbf{x}, t) = f(\mathbf{x}, t), 
		& (\mathbf{x}, t) \in \Omega \times (0, T], \\
		\mathcal{I}u(\mathbf{x}, 0) = u_0(\mathbf{x}), 
		& \mathbf{x} \in \Omega, \\
		\mathcal{B}u(\mathbf{x}, t) = g(\mathbf{x}, t), 
		& (\mathbf{x}, t) \in \partial\Omega \times (0, T],
	\end{cases}
	\label{eq:pde}
\end{equation}
where $\Omega \subset \mathbb{R}^d$ is a bounded spatial domain with boundary $\partial\Omega$, $T > 0$ is the final time.
{Here,} $\mathcal{L}$ is the spatiotemporal differential operator, $\mathcal{I}$ is the initial operator on $\Omega$, $\mathcal{B}$ is the boundary operator on $\partial\Omega$, ${f}$ is the source term, ${u_0}$ is the initial condition, and ${g}$ is the {time-dependent} boundary condition.

Physics-informed neural networks (PINNs)~\cite{MR3881695} provide a mesh-free technique for solving PDEs by embedding the 
physical laws directly into the training objective of a neural network. In this approach, a neural network $u_{\boldsymbol{\theta}_p}(\mathbf{x}, t)$ with trainable parameters $\boldsymbol{\theta}_p$ is employed to approximate the solution $u(\mathbf{x}, t)$ of \eqref{eq:pde} over the entire spatiotemporal domain ${\overline{\Omega}} \times [0, T]$. The network is trained by minimizing the following loss functional:
\begin{equation}
	\mathcal{J}(\boldsymbol{\theta}_p) = \mathcal{J}_r(\boldsymbol{\theta}_p) 
	+ \gamma_{bd}\, \mathcal{J}_{bd}(\boldsymbol{\theta}_p) 
	+ \gamma_{ic}\, \mathcal{J}_{ic}(\boldsymbol{\theta}_p),
	\label{eq:loss}
\end{equation}
where $\gamma_{bd} > 0$ and $\gamma_{ic} > 0$ are penalty parameters, and
\begin{align}
	\mathcal{J}_r(\boldsymbol{\theta}_p) &= \int_0^T \int_{\Omega} 
	|\mathcal{L}u_{\boldsymbol{\theta}_p}(\mathbf{x},t) - 
	f(\mathbf{x},t)|^2 \, d\mathbf{x}\, dt, \\
	\mathcal{J}_{bd}(\boldsymbol{\theta}_p) &= \int_0^T 
	\int_{\partial\Omega} |\mathcal{B}u_{\boldsymbol{\theta}_p}
	(\mathbf{x},t) - g(\mathbf{x},t)|^2 \, d\mathbf{x}\, dt, \\
	\mathcal{J}_{ic}(\boldsymbol{\theta}_p) &= \int_{\Omega} 
	|\mathcal{I}u_{\boldsymbol{\theta}_p}(\mathbf{x},0) - 
	u_0(\mathbf{x})|^2 \, d\mathbf{x}.
\end{align}
In practice, the integrals in \eqref{eq:loss} are approximated using the Monte Carlo method with three sets of collocation points drawn from a uniform distribution: $\mathcal{S}_r = \{(\mathbf{x}^{(i)}_r, t^{(i)}_r)\}_{i=1}^{{n_r}}$ sampled from $\Omega \times (0, T]$, 
$\mathcal{S}_{bd} = \{(\mathbf{x}^{(i)}_{bd}, t^{(i)}_{bd})\}_{i=1}^{{n_{bd}}}$ sampled from $\partial\Omega \times (0, T]$, and 
$\mathcal{S}_{ic} = \{(\mathbf{x}^{(i)}_{ic}, 0)\}_{i=1}^{{n_{ic}}}$ sampled from $\Omega \times \{0\}$, leading to the empirical loss
\begin{align}
	\mathcal{J}^{{n}}(\boldsymbol{\theta}_p) = &\, 
	\frac{1}{{n_r}}\sum_{i=1}^{{n_r}} 
	|\mathcal{L}u_{\boldsymbol{\theta}_p}(\mathbf{x}^{(i)}_r, 
	t^{(i)}_r) - f(\mathbf{x}^{(i)}_r, t^{(i)}_r)|^2 
	\nonumber\\
	&+ \, \frac{\gamma_{bd}}{{n_{bd}}}\sum_{i=1}^{{n_{bd}}} 
	|\mathcal{B}u_{\boldsymbol{\theta}_p}(\mathbf{x}^{(i)}_{bd}, 
	t^{(i)}_{bd}) - g(\mathbf{x}^{(i)}_{bd}, t^{(i)}_{bd})|^2
	\nonumber\\
&+ \, \frac{\gamma_{ic}}{{n_{ic}}}\sum_{i=1}^{{n_{ic}}} 
	|\mathcal{I}u_{\boldsymbol{\theta}_p}(\mathbf{x}^{(i)}_{ic}, 0) 
	- u_0(\mathbf{x}^{(i)}_{ic})|^2.
	\label{eq:empirical_loss}
\end{align}
This discretization transforms the continuous loss functional into a computable form, where the network solution $u_{\boldsymbol{\theta}_p}(\mathbf{x}, t)$ is constrained to satisfy the governing equations at the sampled collocation points. The resulting optimization problem is then solved using gradient-based methods via automatic differentiation~\cite{MR3800512}.

The accuracy of \eqref{eq:empirical_loss} as an approximation of \eqref{eq:loss} depends critically on the choice of collocation points. It has been shown that the total error of the neural network approximation can be decomposed into an approximation error, determined by the network capacity, and a statistical error, arising from the Monte Carlo discretization~\cite{MR4531552, tang2023adversarial}. 
For problems with strongly localized solutions, uniform sampling leads to large statistical errors, particularly in high dimensions, 
which motivates the use of adaptive sampling strategies. It can be shown~\cite{MR4531552} that the statistical error is minimised 
when collocation points are sampled from the distribution proportional to $\mathcal{R}^2(\mathbf{x}, t; \boldsymbol{\theta}_p)$, where 
$\mathcal{R}(\mathbf{x}, t; \boldsymbol{\theta}_p) = \mathcal{L} u_{\boldsymbol{\theta}_p}(\mathbf{x},t) - f(\mathbf{x},t)$ denotes the 
PDE residual, concentrating points where the PDE is most poorly satisfied.

\section{Deep adaptive sampling for spatiotemporal domains}\label{sec:DAS-PINNs}

We apply the DAS-PINNs framework~\cite{MR4531552} to time-dependent high-dimensional PDEs by treating the full spatiotemporal domain $\Omega \times [0, T]$ as a unified high-dimensional space, with the time variable $t$ handled on equal footing with the spatial {variable} $\mathbf{x}$. It alternates between two coupled steps: training the solution network $u_{\boldsymbol{\theta}_p}$ to minimise the loss functional~{\eqref{eq:empirical_loss}}, and training KRnet~\cite{tang2020} to approximate the residual-induced distribution and generate new collocation points concentrated in high-residual regions. KRnet is a normalising flow that defines an invertible transport map $f_{ \boldsymbol{\theta}_k}: \mathbb{R}^{d+1} \to \mathbb{R}^{d+1}$ based on the Knothe--Rosenblatt rearrangement, through which an explicit probability density model $p_{\boldsymbol{\theta}_k}(\mathbf{x}, t)$ is obtained via the change of variables formula
\begin{equation}
	p_{ \boldsymbol{\theta}_k}(\mathbf{x}, t) = p_Z(f_{ \boldsymbol{\theta}_k}(\mathbf{x}, t)) 
	\left|\det \nabla_{(\mathbf{x},t)} f_{ \boldsymbol{\theta}_k}\right|,
	\label{eq:krnet_pdf}
\end{equation}
where $p_Z$ is a standard Gaussian prior. Samples from $p_{ \boldsymbol{\theta}_k}$ are generated efficiently via the inverse map $(\mathbf{x}, t) = f_{ \boldsymbol{\theta}_k}^{-1}(\boldsymbol{\mathrm{z}})$, where $\boldsymbol{\mathrm{z}} \sim p_Z$. Since KRnet is defined on $\mathbb{R}^{d+1}$ and assigns positive density everywhere, it is not directly consistent with the residual-induced distribution, which is zero outside $\Omega \times [0, T]$. To address this, a bounded support mapping is applied to establish a bijection between the bounded spatiotemporal domain and $\mathbb{R}^{d+1}$, ensuring that the generated collocation points remain 
within $\Omega \times [0, T]$.

{For simplicity, we assume that the spatial domain  can be mapped to a  reference domain  $\Omega= [-1,1]^d$, 
by an  affine transformation.} The time variable $t \in [0,T]$ is {similarly}  rescaled to $\tilde{t} \in [-1,1]$ via
\begin{equation}
	\tilde{t} = \frac{2t}{T} - 1,
	\label{eq:time_map}
\end{equation}
so that the full spatiotemporal domain $\Omega \times [0,T]$ is embedded in $\Omega_n = [-1,1]^{d+1}$.

Let $\mathcal {D}_n = (-(1+\delta), 1+\delta)^{d+1}$ with $\delta > 0$ be a slightly enlarged domain containing $[-1,1]^{d+1}$. To map $\mathcal {D}_n $ to $\mathbb{R}^{d+1}$, the following bijection is applied componentwise to each 
spatiotemporal variable $\xi_i \in (-(1+\delta), 1+\delta)$:
\begin{equation}
	\Phi(\xi_i) = \log\frac{\xi_i + (1+\delta)}{(1+\delta) - \xi_i},
	\label{eq:logit_map}
\end{equation}
which maps $(-(1+\delta), 1+\delta)$ to $\mathbb{R}$ and has positive derivative everywhere. The composition mapping $f_{\boldsymbol{\theta}_k} \circ \Phi$ then defines an explicit density model on $\mathcal {D}_n $ as
\begin{equation}
	\hat{p}_{ \boldsymbol{\theta}_k}(\boldsymbol{\xi}) = p_{ \boldsymbol{\theta}_k}(\Phi(\boldsymbol{\xi})) 
	|\nabla_{\boldsymbol{\xi}} \Phi|,
	\label{eq:bounded_pdf}
\end{equation}
where $\boldsymbol{\xi} = (\xi_1, \ldots, \xi_{d+1}) \in \mathcal{D}_n$ denotes the full normalised spatiotemporal input in $\mathcal{D}_n $. In our experiments, we set $\delta = 0.01$. 

To ensure that $\hat{p}_{ \boldsymbol{\theta}_k}$ and the residual-induced distribution share the same support $\mathcal {D}_n $, a cutoff function $h: \mathcal {D}_n  \to [0,1]$ is introduced following~\cite{MR4531552}, defined component-wise as $h(\boldsymbol{\xi}) = \prod_{i=1}^{d+1} h_\delta(\xi_i)$, where $h_\delta(\xi_i) = 1$ on $[-1,1]$, tapers linearly to zero on 
$\mathcal {D}_n  \setminus [-1,1]^{d+1}$, and vanishes outside $\mathcal {D}_n $. The modified residual distribution is then
\begin{equation}
\widetilde{\mathcal{R}}(\boldsymbol{\xi}) \propto 
\mathcal{R}^2((\mathbf{x},t); \boldsymbol{\theta}_p)\, h(\boldsymbol{\xi}),
	\label{eq:modified_residual}
\end{equation}
where $(\mathbf{x}, t) \in \Omega \times [0,T]$ {denotes the corresponding spatiotemporal point}. 
KRnet is trained by solving
\begin{equation}
	\boldsymbol{\theta}_k^* = \arg\min_{\boldsymbol{\theta}_k}\, 
	\mathds{D}_{\mathrm{KL}}\!\left(\widetilde{\mathcal{R}}(\boldsymbol{\xi})\, 
	\Big\|\, \hat{p}_{\boldsymbol{\theta}_k}(\boldsymbol{\xi})\right),
	\label{eq:kl_opt}
\end{equation}
where $\mathds{D}_{\mathrm{KL}}(\cdot\|\cdot)$ denotes the {Kullback--Leibler (KL)}  divergence. Since minimising the KL divergence is equivalent to minimising the cross entropy between $\widetilde{\mathcal{R}}$ and $\hat{p}_{\boldsymbol{\theta}_k}$, the cross entropy is approximated via importance sampling as
\begin{equation}
	\mathcal{H}(\widetilde{\mathcal{R}}, \hat{p}_{\boldsymbol{\theta}_k}) 
	\approx -\frac{1}{{n_r}}\sum_{i=1}^{{n_r}} 
	\frac{\widetilde{\mathcal{R}}(\boldsymbol{\xi}_{\mathcal {D}_n}^{i})}
	{\hat{p}_{\tilde{\boldsymbol{\theta}}_k}(\boldsymbol{\xi}_{\mathcal {D}_n}^{i})}\,
	\log \hat{p}_{\boldsymbol{\theta}_k}(\boldsymbol{\xi}_{\mathcal {D}_n}^{i}),
	\label{eq:cross_entropy_approx}
\end{equation}
where $\hat{p}_{\tilde{\boldsymbol{\theta}}_k}$ is a density model with known parameters $\tilde{\boldsymbol{\theta}}_k$, and its samples $\{\boldsymbol{\xi}_{\mathcal {D}_n}^{i}\}_{i=1}^{{n_r}}$ are generated as
\begin{equation}
\boldsymbol{\xi}_{\mathcal {D}_n}^{i} = \Phi^{-1} \circ 
	f_{\tilde{\boldsymbol{\theta}}_k}^{-1}(\mathbf{z}^{(i)}), 
	\quad \mathbf{z}^{(i)} \sim p_Z,
	\label{eq:sample_gen}
\end{equation}
with $\mathbf{z}^{(i)}$ sampled from the Gaussian prior. How $\tilde{\boldsymbol{\theta}}_k$ is chosen at each stage will become clear in {the next section.}

\subsection{{Adaptive sampling strategy}}
\label{sec:asp}

The key adaptation to the time-dependent setting is that three separate sets of collocation points are maintained throughout training:
\begin{align}
	\mathcal{S}_r^j &= \{(\mathbf{x}_r^{(i)}, t_r^{(i)})\}_{i=1}^{{n_r}^j}, \quad
	\mathcal{S}_{bd}^j = \{(\mathbf{x}_{bd}^{(i)}, t_{bd}^{(i)})\}_{i=1}^{{n_{bd}}^j},  \quad
	\mathcal{S}_{ic}^j = \{(\mathbf{x}_{ic}^{(i)},0)\}_{i=1}^{{n_{ic}}^j}, 
\end{align}
where $j = 0,1, \ldots, {n}_{\mathrm{stage}}-1$ denotes the adaptive stage index. At stage $j = 0$, all three sets are initialised with uniform samples from their respective domains: $\mathcal{S}_r^0$ from $\Omega \times (0, T]$, $\mathcal{S}_{bd}^0$ from $\partial\Omega \times (0, T]$, and $\mathcal{S}_{ic}^0$ from $\Omega \times \{0\}$. The solution network is trained by minimising the empirical loss~\eqref{eq:empirical_loss} to obtain $u_{\boldsymbol{\theta}_p^{*,(1)}}(\mathbf{x}, t)$. With $u_{\boldsymbol{\theta}_p^{*,(1)}}(\mathbf{x}, t)$ in hand, KRnet is trained on $\mathcal{S}_r^0$ by minimising cross entropy~\eqref{eq:cross_entropy_approx} using uniform samples to get $\hat{p}_{\boldsymbol{\theta}_k^{*,(1)}}(\boldsymbol{\xi})$. KRnet generates $n_r$ new candidate points $\{\boldsymbol{\xi}_{\mathcal{D}_n}^{i,(1)}\}_{i=1}^{n_r}$ in $\mathcal{D}_n$ via~\eqref{eq:sample_gen} using the optimal parameters $\boldsymbol{\theta}_k^{*,(1)}$. These points are processed through a projection and classification step described later and then denormalised to 
$\mathcal{D}_n^d = (-(1+\delta), 1+\delta)^d \times \left(-\frac{\delta T}{2}, T + \frac{\delta T}{2}\right)$ to obtain the new set $\mathcal{S}_r^{a,1} = \{(\mathbf{x}_r^{(i)}, t_r^{(i)})\}_{i=1}^{n_r}$. The set is then updated as $\mathcal{S}_r^1 = \mathcal{S}_r^{0} \cup \mathcal{S}_r^{a,1}$, and the solution network is updated to obtain $u_{\boldsymbol{\theta}_p^{*,(2)}}(\mathbf{x}, t)$.

In general, at each subsequent stage $j \geq 1$, we use 
\begin{align}
	\mathcal{S}_r^j &= \mathcal{S}_r^{j-1} \cup \mathcal{S}_r^{a,j}, \quad
	\mathcal{S}_{bd}^j = \mathcal{S}_{bd}^{j-1} \cup \mathcal{S}_{bd}^{a,j}, \quad
	\mathcal{S}_{ic}^j = \mathcal{S}_{ic}^{j-1} \cup \mathcal{S}_{ic}^{a,j},
\end{align}
to obtain $u_{\boldsymbol{\theta}_p^{*,(j+1)}}(\mathbf{x}, t)$ by solving
\begin{equation}
	\boldsymbol{\theta}_p^{*,(j+1)} = \arg\min_{\boldsymbol{\theta}_p}\, 
	\mathcal{J}^{{n}}(\boldsymbol{\theta}_p; \mathcal{S}_r^j, 
	\mathcal{S}_{bd}^j, \mathcal{S}_{ic}^j),
	\label{eq:pde_opt_j}
\end{equation}
where $\boldsymbol{\theta}_p$ is initialised as $\boldsymbol{\theta}_p^{*,(j)}$ and $\mathcal{J}^{{n}}$ is given in~\eqref{eq:empirical_loss}.

Further, with $\hat{p}_{\boldsymbol{\theta}_k^{*,(j)}}(\boldsymbol{\xi})$, KRnet is updated via 
\begin{equation}
	\boldsymbol{\theta}_k^{*,(j+1)} = \arg\min_{\boldsymbol{\theta}_k}\,  -\frac{1}{{n_r}^j}\sum_{i=1}^{{n_r}^j} 
	\frac{\mathcal{R}^2\left((\mathbf{x}_{{\mathcal {D}_n^d}}^{i, (j)}, t_{{\mathcal {D}_n^d}}^{i, (j)}); 	\boldsymbol{\theta}_p^{*,(j+1)}\right)h\left(\boldsymbol{\xi}_{\mathcal {D}_n}^{i, (j)}\right)}
	{\hat{p}_{\boldsymbol{\theta}_k^{*,(j)}}(\boldsymbol{\xi}_{\mathcal {D}_n}^{i, (j)})}\,
	\log \hat{p}_{\boldsymbol{\theta}_k}(\boldsymbol{\xi}_{\mathcal {D}_n}^{i, (j)}),
	\label{eq:cross_entropy_approx1}
\end{equation}
where the parameters $\tilde{\boldsymbol{\theta}}_k$ in~\eqref{eq:cross_entropy_approx} are set to $\boldsymbol{\theta}_k^{*,(j)}$, the optimal KRnet parameters from the previous stage. 

\subsection{{Sampling point projection strategy}}
\label{sec:spp}
A new set of points $\{\boldsymbol{\xi}_{\mathcal{D}_n}^{i,(j+1)}\}_{i=1}^{n_r}$ is then generated via~\eqref{eq:sample_gen} using the optimal parameters $\boldsymbol{\theta}_k^{*,(j+1)}$. These points may lie anywhere in $\mathcal{D}_n$, including outside $\Omega_n = [-1,1]^{d+1}$. They are processed through a projection and classification step described below and then denormalised to $\mathcal{D}_n^d$ to obtain the new sets $\mathcal{S}_r^{a,j+1}$, $\mathcal{S}_{bd}^{a,j+1}$ and $\mathcal{S}_{ic}^{a,j+1}$. 

Points in $\mathcal{D}_n \setminus \Omega_n$ are projected onto the boundary of $\Omega_n$ via an entry-wise clipping operator 
following~\cite{MR4531552}. The key adaptation to the time-dependent setting is the subsequent classification of projected points: since the last coordinate $\xi_{d+1}$ corresponds to the normalised time $\tilde{t}$, a projected point with $\xi_{d+1} = -1$ corresponds to $t = 0$ and is assigned to $\mathcal{S}_{ic}^{a,j+1}$; a point with any spatial coordinate $|\xi_i| = 1$ for $i = 1, \ldots,d$ and $\xi_{d+1} \neq -1$ is assigned to $\mathcal{S}_{bd}^{a,j+1}$; remaining projected points are assigned back to $\mathcal{S}_r^{a,j+1}$. If the projected points are insufficient to fill $\mathcal{S}_{bd}^{a,j+1}$ or $\mathcal{S}_{ic}^{a,j+1}$ to size $n_r$, the remainder is filled with fresh uniform samples from $\partial\Omega \times (0,T]$ and $\Omega \times \{0\}$ respectively. This procedure is repeated until the stopping criterion is satisfied or the maximum number of 
stages ${n}_{\mathrm{stage}}$ is reached.

\section{Results and discussion}\label{sec:Results}

We present numerical experiments on a range of time-dependent PDEs to demonstrate the effectiveness of the DAS-PINNs approach extended to the spatiotemporal domain. The experiments are organized into two groups: low-dimensional problems in $\mathbb{R}^2 \times [0,T]$ and high-dimensional problems in $\mathbb{R}^d \times [0,T]$ with $d$ up to $8$. In all experiments, the density model 
$p_{\boldsymbol{\theta}_k}$ is implemented using KRnet~{\cite{tang2020}}, and all networks are trained using the Adam optimizer~\cite{kingma2014adam} with a fixed learning rate of $10^{-4}$. The solution network architecture and training configuration are specified for each problem individually. The accuracy of the approximate solution is measured using three {pointwise errror} metrics: the relative $L_2$ error ($\mathcal{E}_2$), the $L_{\infty}$ error ($\mathcal{E}_\infty$), 
and the mean square error ($\mathcal{E}_{\mathrm{mse}}$), defined respectively as
\begin{equation}
	\mathcal{E}_2 = \frac{\|u_{\boldsymbol{\theta}_p} - u\|_2}{\|u\|_2}, 
	\qquad
	\mathcal{E}_\infty = \|u_{\boldsymbol{\theta}_p} - u\|_\infty,
	\qquad
	\mathcal{E}_{\mathrm{mse}}= \frac{1}{{n}_v}\sum_{i=1}^{{n}_v}
	|u_{\boldsymbol{\theta}_p}(\mathbf{x}^{(i)}, t^{(i)}) - 
	u(\mathbf{x}^{(i)}, t^{(i)})|^2,
\end{equation}
where $u$ denotes the exact solution and ${n}_v$ is the number of validation points {in  $\mathbb{R}^d \times [0,T]$}.

\subsection{Low-dimensional problems}
In this section, we consider three time-dependent {PDE} problems in $\mathbb{R}^2 \times [0, T]$: a moving transient Gaussian peak, a rotation equation, and {the}  viscous Burgers' equation.

\subsubsection{{Transient Gaussian peak}}\label{sec:tgp}
We consider the following parabolic {PDE} problem:
\begin{equation}
	\begin{cases}
		\dfrac{\partial u}{\partial t} - \Delta u = f(\mathbf{x}, t), 
		& (\mathbf{x}, t) \in \Omega \times (0, T], \\
		u(\mathbf{x}, 0) = 0, 
		& \mathbf{x} \in \Omega, \\
		u(\mathbf{x}, t) = g(\mathbf{x}, t), 
		& (\mathbf{x}, t) \in \partial\Omega \times (0, T],
	\end{cases}
	\label{eq:moving_transient}
\end{equation}
where $\Omega = [-1,1]^2$, $T = 1$. The source term $f$ is chosen such 
that the exact solution is:
\begin{equation}
	u(\mathbf{x}, t) = \exp\bigl(-\alpha\{(x_1 + 0.5 - 0.5t)^2 
	+ x_2^2\}\bigr) \cdot (1 - e^{-5t}),
	\label{eq:transient_peak}
\end{equation}
with $\alpha = 100$. The peak centre moves linearly from $(-0.5, 0)$ at $t = 0$ to $(0, 0)$ at $t = 1$ with velocity $c = 0.5$ in the $x_1$ direction, while the solution {evolves} from zero at $t=0$ toward a steady-state Gaussian profile.

\begin{figure}[tbhp!]
	\centering
	\includegraphics[width=0.95\linewidth]{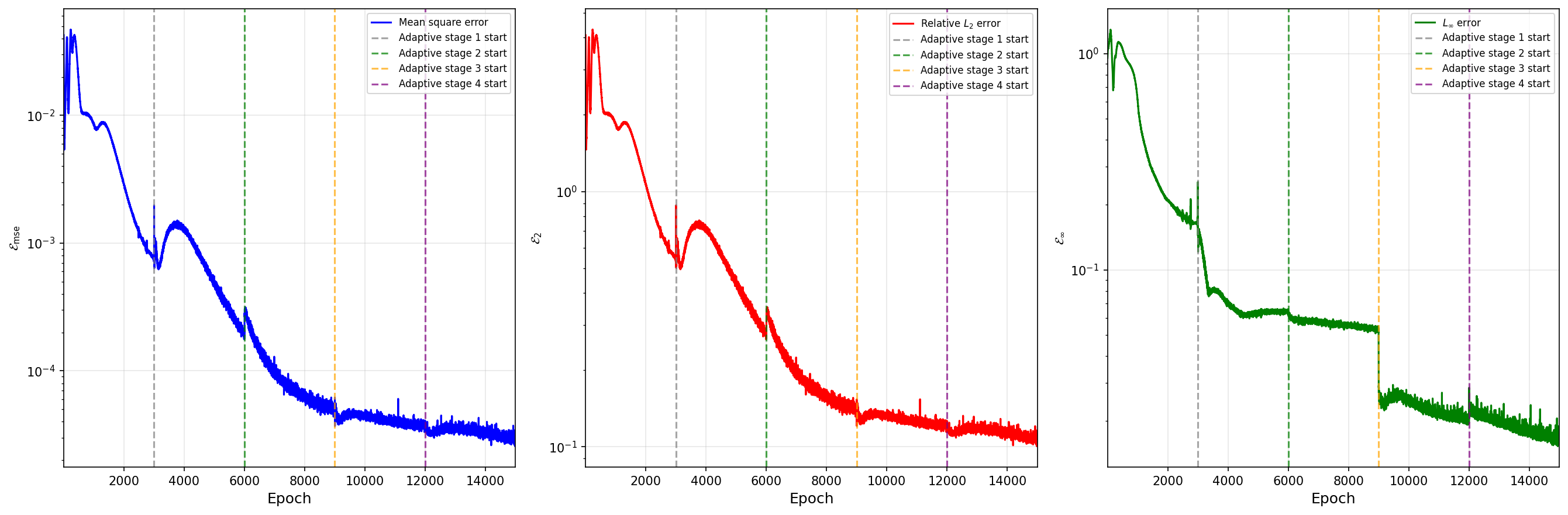}
	\caption{{Moving Gaussian peak problem \eqref{eq:moving_transient}}: evolution 
		of the mean square error ($\mathcal{E}_{mse}$), relative $L_2$ error 
		($\mathcal{E}_2$), and $L_{\infty}$ error ($\mathcal{E}_\infty$) per epoch. 
		The initial stage is followed by four adaptive stages. Dashed vertical 
		lines mark the start of each adaptive stage.}\label{fig:prob1_error_vs_epoch}
\end{figure}

\begin{figure}[tbhp!]
	\centering
	\includegraphics[width=0.95\linewidth]{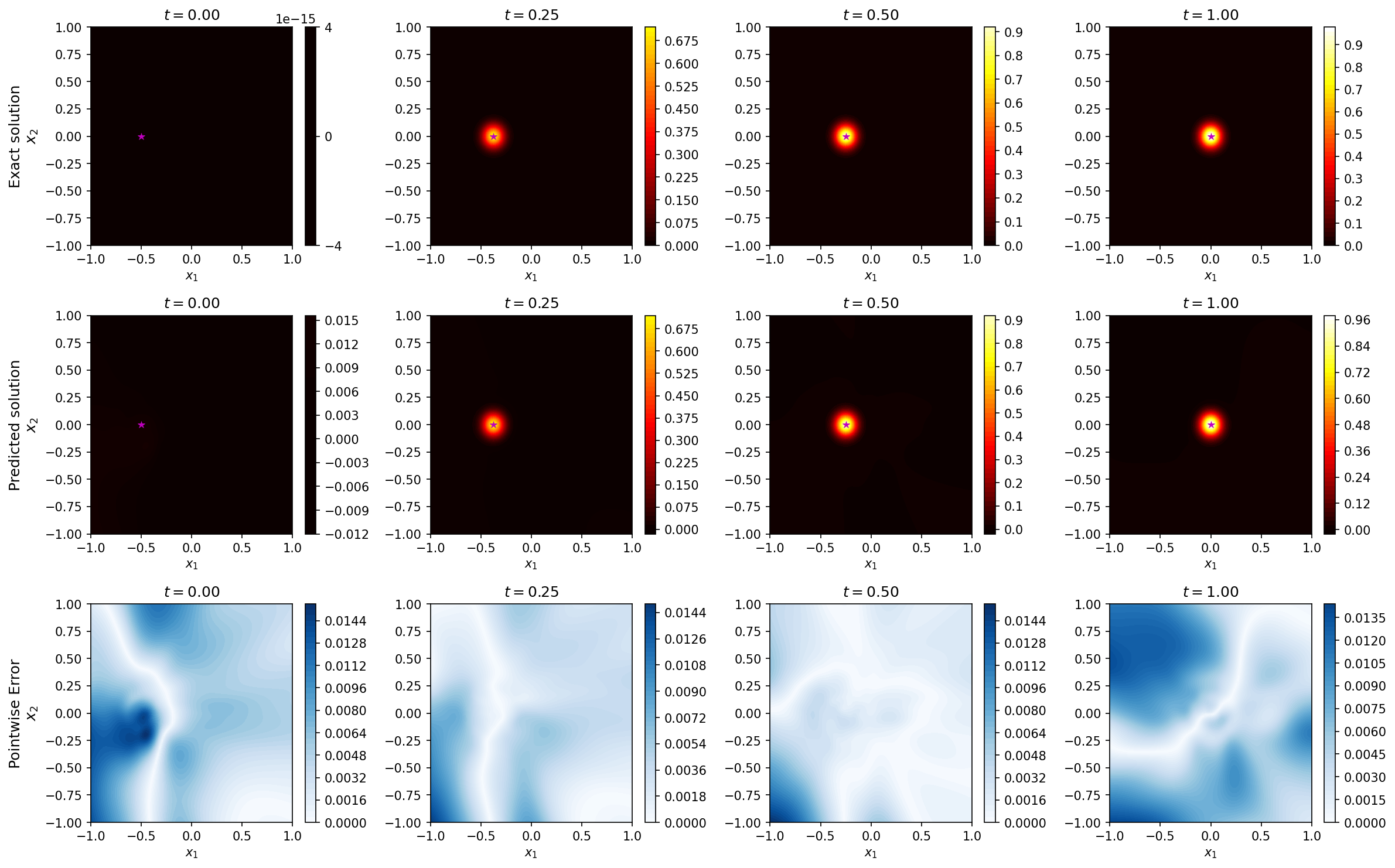}
	\caption{ {Moving Gaussian peak problem \eqref{eq:moving_transient}}: exact solution (top row), 
		predicted solution (middle row), and pointwise error (bottom row) 
		at {four  time levels $t \in \{0, 0.25, 0.50, 1.00\}$}.}\label{fig:prob1_solution_slices}
\end{figure}

The solution network and KRnet are trained with $n_{\text{train}} = 2000$ collocation points and batch size $500$, using a learning rate of $10^{-4}$. The solution network is trained for $3000$ epochs per stage and the flow network for $3000$ epochs per stage, over one initial stage followed by four adaptive stages, giving a total of $180{,}000$ training iterations. The penalty parameters are set to $ \gamma_{bd} = 1.0$ and $\gamma_{ic} = 1.0$. The solution network uses $6$ fully connected layers with $32$  hidden neurons per layer. For KRnet, $6$ affine coupling layers are used, each consisting of two fully connected layers with $24$ neurons per layer. The model is evaluated on a tensor product validation grid of $100 \times 100$ spatial points over $[-1,1]^2$ at $5$ uniformly spaced time levels $t \in \{0, 0.25, 0.5, 0.75, 1.0\}$, giving a total of $50{,}000$ validation points.

Figure~\ref{fig:prob1_error_vs_epoch} shows the evolution of ${\mathcal{E}_{\mathrm{mse}}}$, $\mathcal{E}_2$, and $\mathcal{E}_\infty$ across all training epochs. All three metrics decrease consistently across the five stages, with each adaptive stage producing a  notable reduction in error. The final minimum values achieved are ${\mathcal{E}_{\mathrm{mse}}} = 2.11 \times 10^{-5}$, $\mathcal{E}_2 = 9.14 \times 10^{-2}$, and $\mathcal{E}_\infty = 1.32 \times 10^{-2}$. Figure~\ref{fig:prob1_solution_slices} shows the exact solution, predicted  solution, and pointwise error at  {four} time levels. At $t = 0$ the solution is identically zero due to the factor $(1 - e^{-5t})$, with the star marking the initial peak position from which the Gaussian profile emerges for $t > 0$. The predicted solution closely matches the exact solution across all time levels, with the pointwise error remaining uniformly small throughout the evolution. This problem combines temporal transience with spatial motion, and the results demonstrate that the DAS-PINNs handles both aspects effectively.

\begin{figure}[tbhp!]
	\centering
	\includegraphics[width=0.85\linewidth]{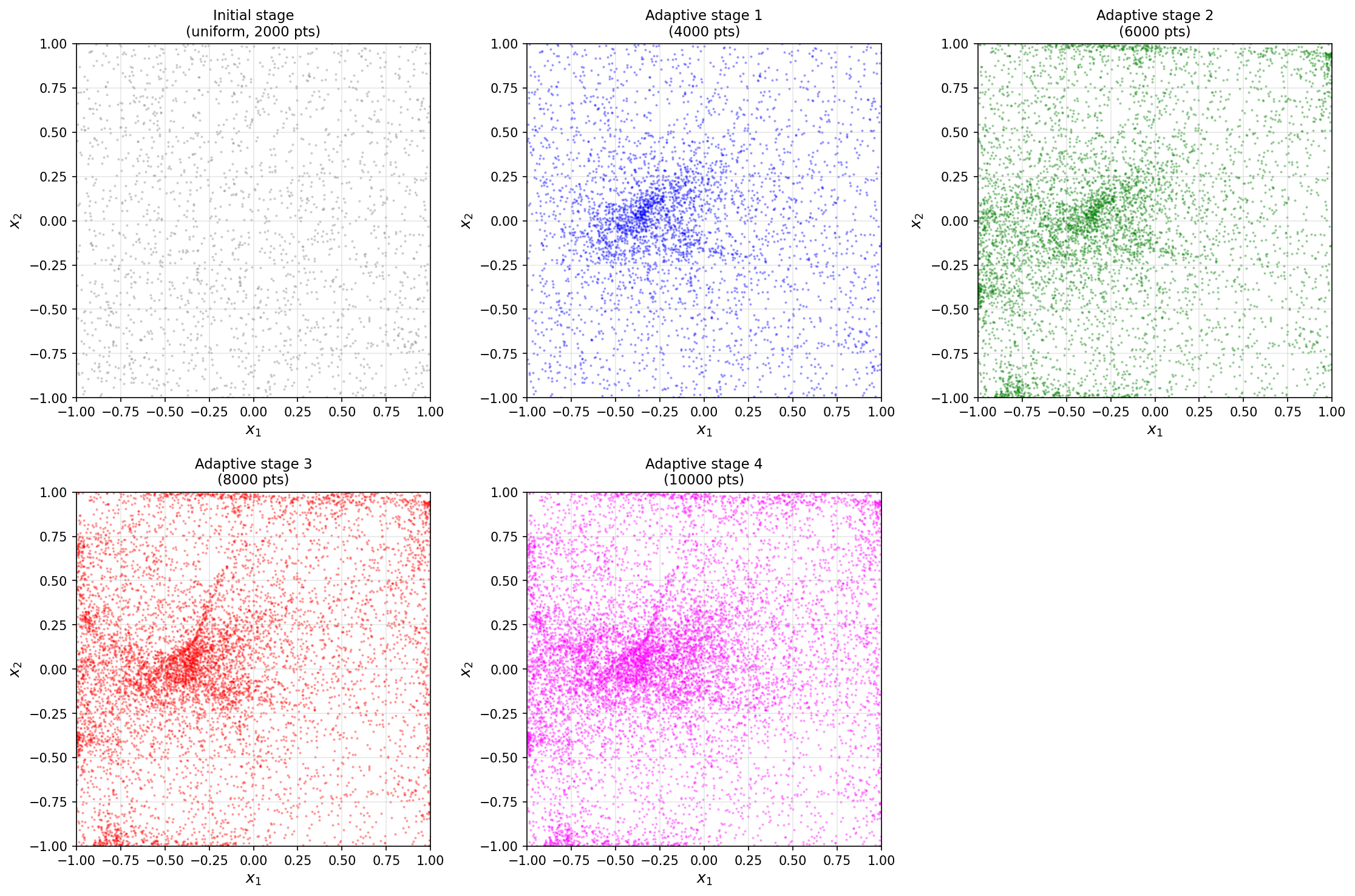}
	\caption{{Moving Gaussian peak problem \eqref{eq:moving_transient}}: spatial distribution of 
		collocation points at each stage of DAS-PINNs. The training set 
		grows from $2000$ points in the initial stage to $10{,}000$ points 
		by the final adaptive stage.}\label{fig:prob1_spatial_points}
\end{figure}

\begin{figure}[tbhp!]
	\centering
	\includegraphics[width=0.85\linewidth]{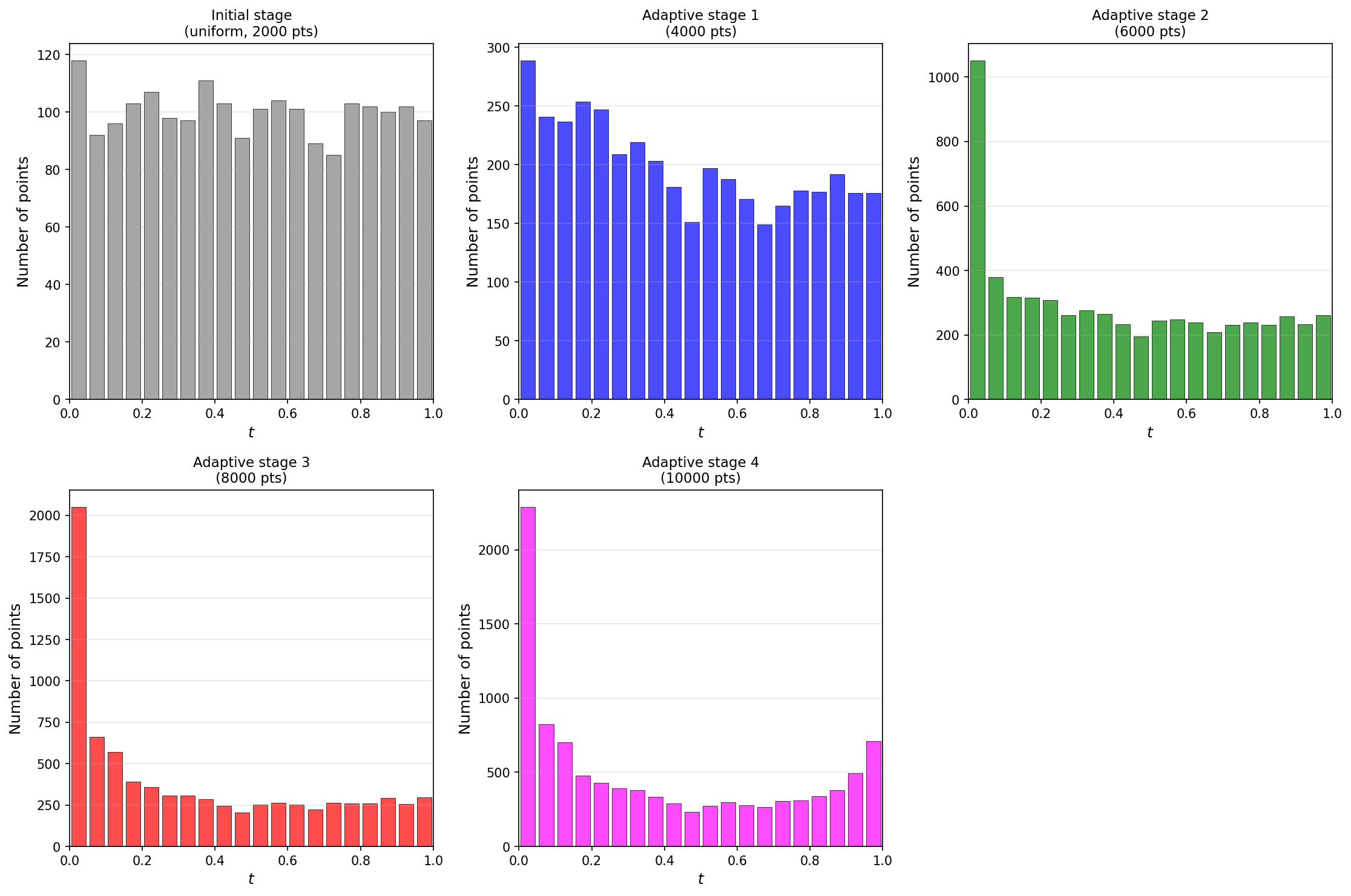}
	\caption{{Moving Gaussian peak problem \eqref{eq:moving_transient}}:  temporal histogram of 
		collocation points at each stage of DAS-PINNs.}\label{fig:prob1_time_points}
\end{figure}

Figure~\ref{fig:prob1_spatial_points} shows the spatial distribution of collocation points across all five stages. Starting from a uniform distribution in the initial stage, DAS-PINNs progressively concentrates points near the peak trajectory, reflecting the residual-driven nature of the adaptive sampling. As the training advances, the point distribution spreads further across the domain, indicating that the residual is becoming more uniform---consistent with the error reduction observed in Figure~\ref{fig:prob1_error_vs_epoch}. Figure~\ref{fig:prob1_time_points} shows the temporal histogram of interior collocation points at each stage. In the initial stage, points are distributed uniformly across the time interval $(0, 1]$. As training progresses, DAS-PINNs concentrates more interior points near $t = 0$. This is physically meaningful, the factor $(1 - e^{-5t})$ in the exact solution causes the solution to grow rapidly from zero at early times, producing large temporal gradients and consequently large PDE residuals near $t = 0$. KRnet naturally identifies these high-residual regions and concentrates interior points there without any explicit guidance.

\subsubsection{{A rotating peak problem}}\label{sec:rotation}
We consider the following pure transport {problem} in which a Gaussian peak rotates about the origin:
\begin{equation}
	\begin{cases}
		\dfrac{\partial u}{\partial t} = \sin t \, \dfrac{\partial u}{\partial x_1} 
		- \cos t \, \dfrac{\partial u}{\partial x_2}, 
		& (\mathbf{x}, t) \in \Omega \times (0, T], \\[6pt]
		u(\mathbf{x}, 0) = \exp\!\left(-100
		\left[(x_1-1)^2 + x_2^2\right]\right), 
		& \mathbf{x} \in \Omega, \\[6pt]
		u(\mathbf{x}, t) = \exp\!\left(-100
		\left[(x_1-\cos t)^2 + (x_2 - \sin t)^2\right]\right), 
		& (\mathbf{x}, t) \in \partial\Omega \times (0, T],
	\end{cases}
	\label{eq:rotation}
\end{equation}
where $\Omega = [-1,1]^2$ and $T = 1.0$. The exact solution is
\begin{equation}
u(\mathbf{x}, t) = \exp\!\left(-100\left[(x_1 - \cos t)^2 + (x_2 - \sin t)^2\right]\right).
\end{equation}

\begin{figure}[tbhp!]
\centering
\includegraphics[width=0.95\linewidth]{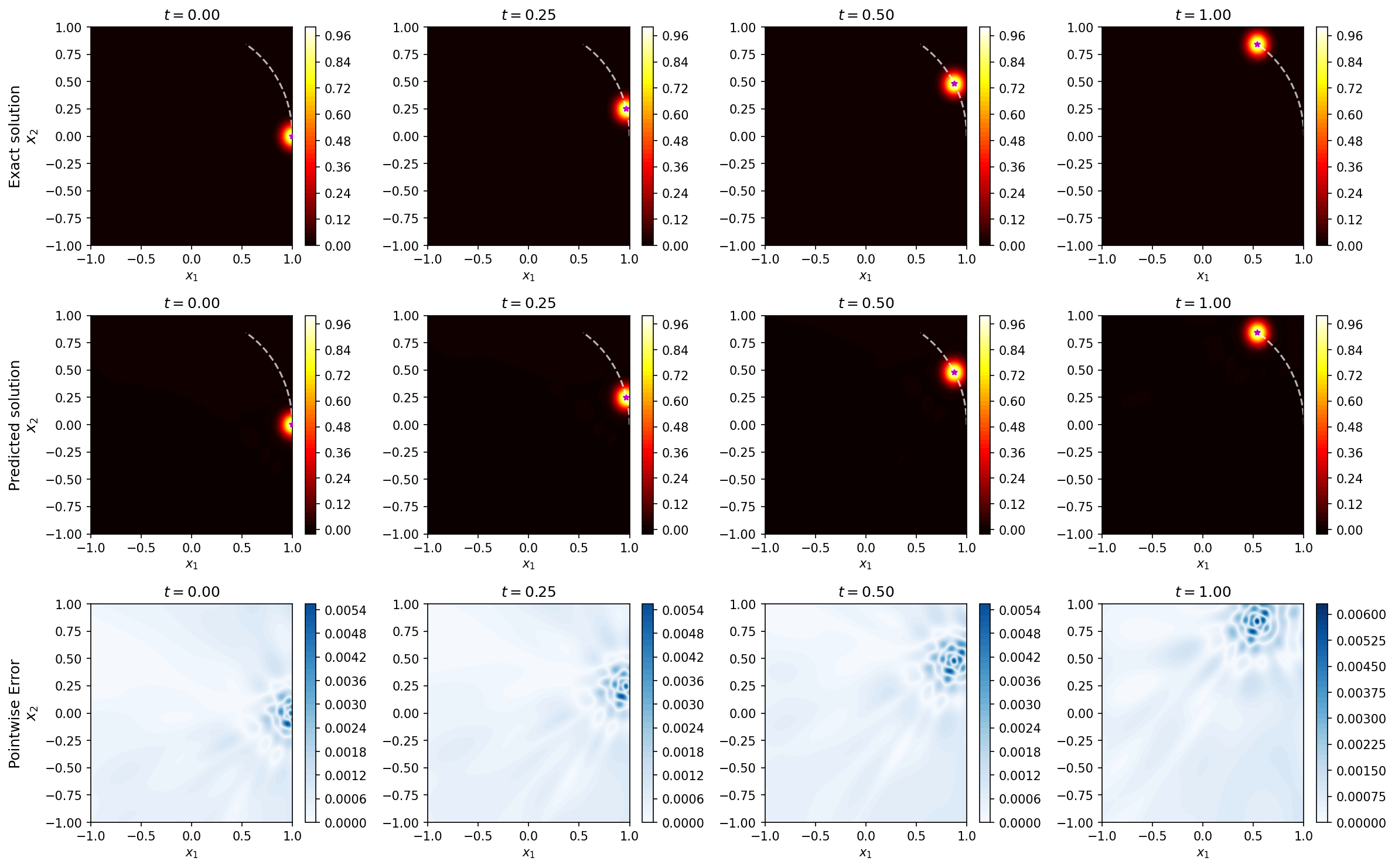}
\caption{{Rotation problem \eqref{eq:rotation}}: exact solution (top row), predicted solution (middle row), and pointwise error (bottom row) at {four time levels $t \in \{0, 0.25, 0.50, 1.00\}$}.}\label{fig:prob2_solution_slices}
\end{figure}

\begin{table}[tbhp!]
\centering
\caption{{Rotation problem \eqref{eq:rotation}}: minimum $\mathcal{E}_2$ and $\mathcal{E}_\infty$ at each adaptive stage.}\label{tab:prob20}
\begin{tabular}{lcccc}
\toprule
& \textbf{initial stage} & \textbf{adaptive stage 1} 
& \textbf{adaptive stage 2} & \textbf{adaptive stage 3} \\
\midrule
Iter range & $1$--$12000$ & $12001$--$36000$ & $36001$--$72000$ & $72001$--$120000$ \\
$\mathcal{E}_2$ min & $1.71\times10^{-1}$ & $2.32\times10^{-2}$ & $1.20\times10^{-2}$ & $8.65\times10^{-3}$ \\
$\mathcal{E}_\infty$ min & $3.17\times10^{-1}$ & $1.67\times10^{-2}$ & $8.38\times10^{-3}$ & $6.59\times10^{-3}$ \\
\bottomrule
\end{tabular}
\end{table}

The peak centre follows the circular trajectory $(\cos t, \sin t)$, starting at $(1, 0)$ at $t = 0$ and tracing a partial arc by $t = 1$.  This equation presents a distinctive challenge for vanilla PINNs, as it is a pure transport problem with no diffusion or source term, meaning there is no natural smoothing mechanism to aid the network in learning the solution. The Gaussian peak is highly localised with coefficient ${100}$ in the exponent, and follows a circular trajectory rather than a simple linear path, making it more demanding for the adaptive sampling to track.
	
The training follows the same configuration as the previous problem, except that the solution network uses $64$ hidden neurons per layer and training runs over one initial stage followed by three adaptive stages, giving a total of $120{,}000$ training iterations. The model is evaluated on a tensor product validation grid of $200 \times 200$ spatial points over $[-1,1]^2$ at $5$ uniformly spaced time levels, giving a total of $200{,}000$ validation points.

Table~\ref{tab:prob20} summarises the minimum $\mathcal{E}_2$ and $\mathcal{E}_\infty$ achieved at each stage. Both errors decrease consistently as the adaptive stages progress, with the most significant reduction occurring after the first adaptive stage. By the final stage, $\mathcal{E}_2$ and $\mathcal{E}_\infty$ reach $8.65 \times 10^{-3}$ and $6.59 \times 10^{-3}$ respectively, demonstrating that DAS-PINNs effectively tracks the rotating peak throughout the time interval. 

\begin{figure}[tbhp!]
\centering
\includegraphics[width=0.65\linewidth]{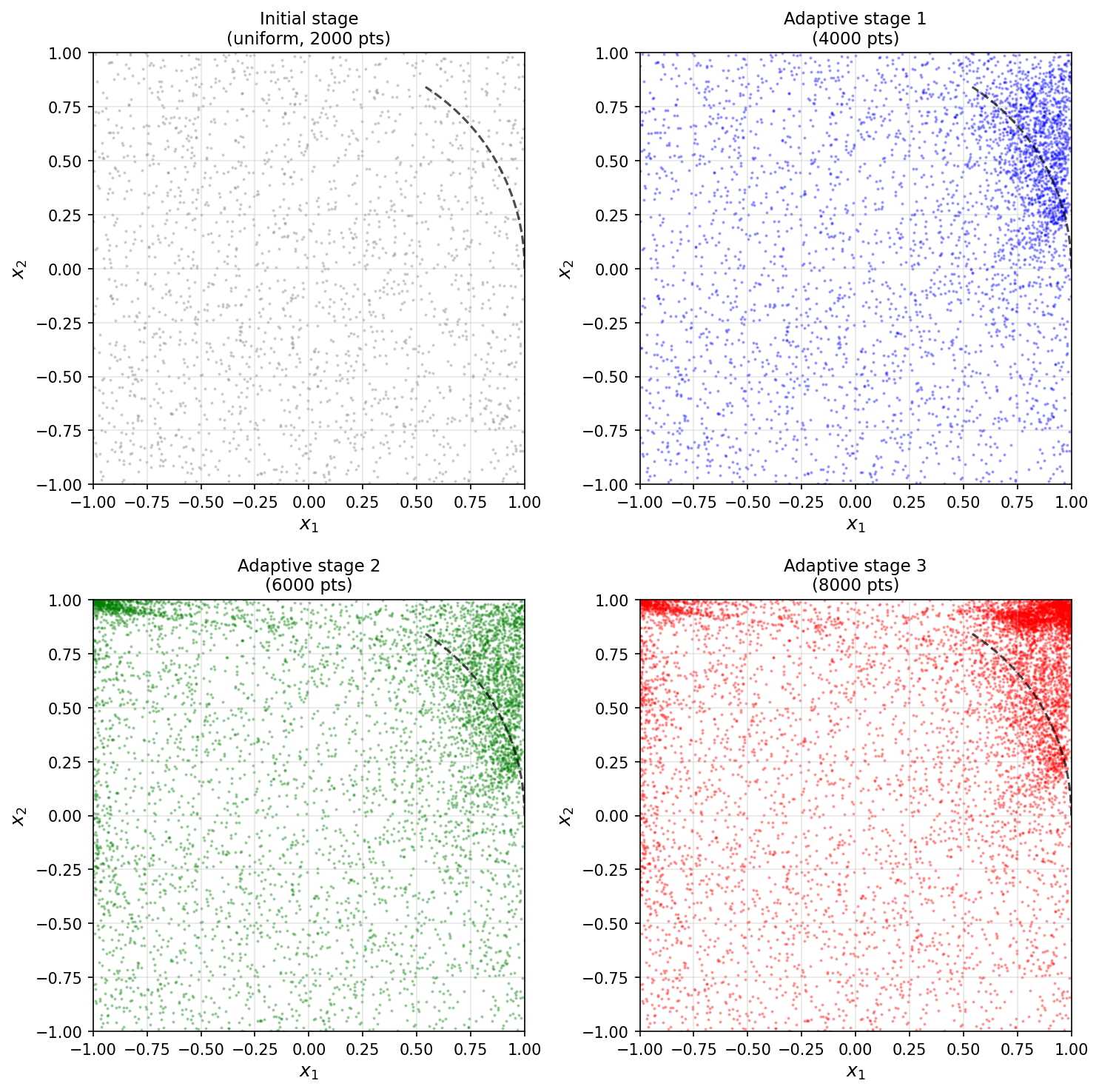}
\caption{{Rotation problem \eqref{eq:rotation}}: spatial distribution of collocation points at each stage of DAS-PINNs. The dashed line shows the circular trajectory of the peak centre.}\label{fig:prob2_spatial_points}
\end{figure}

Figure~\ref{fig:prob2_solution_slices} shows the exact solution, predicted solution, and pointwise error at {four} time levels. The dashed line shows the circular trajectory of the peak centre, and the star marker indicates the peak location at each time level. The network accurately captures the rotating Gaussian peak throughout its circular trajectory, with the pointwise error remaining small across all time levels. Figure~\ref{fig:prob2_spatial_points} shows the spatial distribution of collocation points across all stages. Starting from a uniform distribution in the initial stage, DAS-PINNs progressively concentrates points along the circular trajectory of the rotating peak, reflecting the residual-driven nature of the adaptive sampling. As training advances, the point distribution spreads further across the domain, indicating that the residual is becoming more uniform. 

\begin{figure}[tbhp!]
\centering
\includegraphics[width=0.65\linewidth]{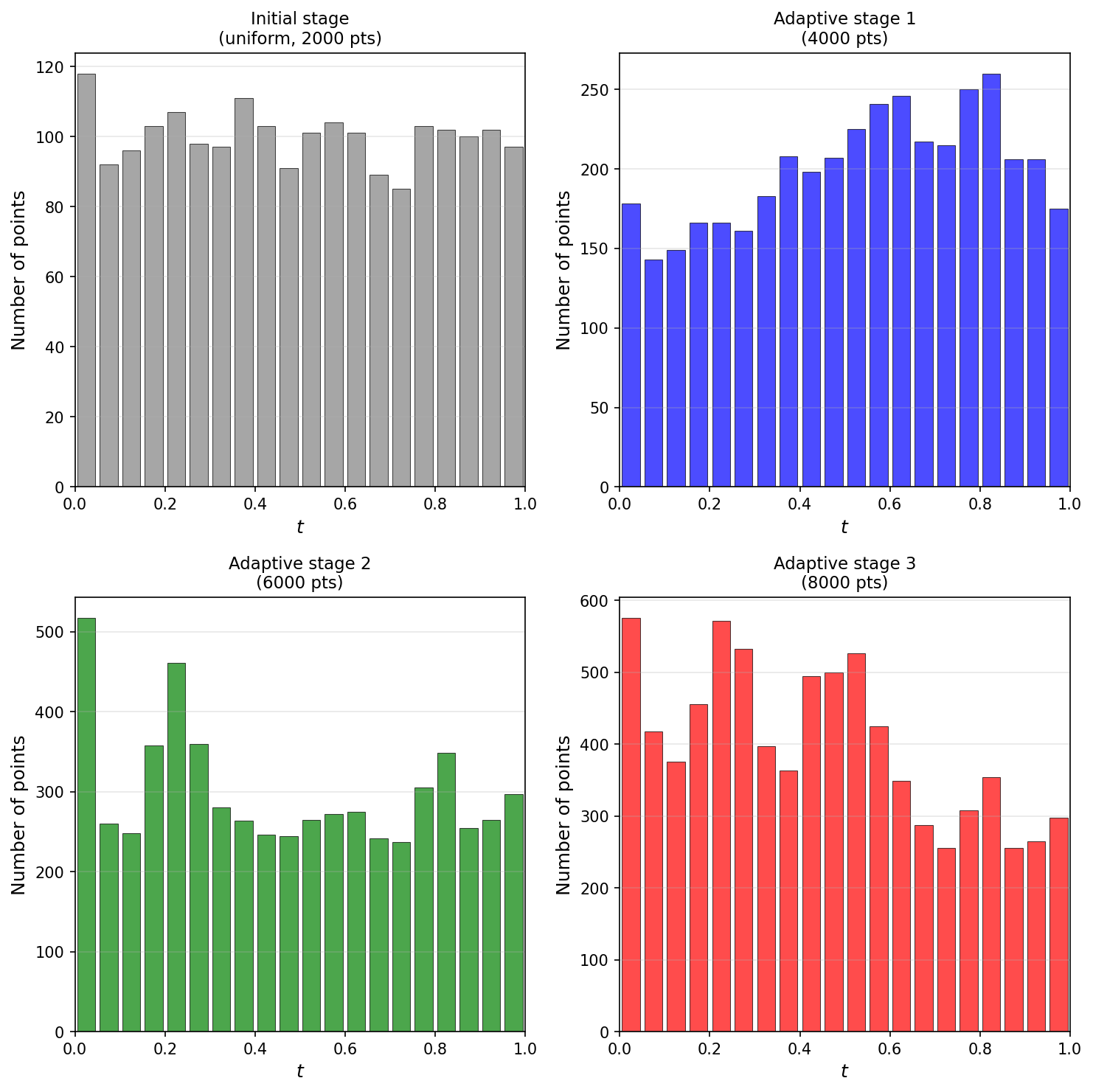}
\caption{{Rotation problem \eqref{eq:rotation}}: temporal histogram of newly added collocation 
points at each adaptive stage of DAS-PINNs.}\label{fig:prob2_time_points}
\end{figure}

Figure~\ref{fig:prob2_time_points} shows the temporal histogram of newly added interior collocation points at each adaptive stage. In the initial stage, points are distributed uniformly across time. As training progresses, the temporal concentration shifts across stages, reflecting the fact that as the residual is reduced in one temporal region, KRnet shifts its attention to the next high-residual region. This behaviour is consistent with the rotating nature of the peak as there is no transient growth factor to bias the sampling toward early times, so KRnet distributes points more evenly across the time interval in response to the evolving residual profile.
	
\subsubsection{{A nonlinear problem}}\label{sec:Burger}

We consider the following {nonlinear problem (viscous Burgers' equation in $\mathbb{R}^2$)}:
\begin{equation}
\begin{cases}
\dfrac{\partial u}{\partial t} = 0.001\!\left(\dfrac{\partial^2 u}{\partial x_1^2} + \dfrac{\partial^2 u}{\partial x_2^2}\right) - u\!\left(\dfrac{\partial u}{\partial x_1} + \dfrac{\partial u}{\partial x_2}\right), & (\mathbf{x}, t) \in \Omega \times (0, T], \\[8pt]
u(\mathbf{x}, 0) = \dfrac{1}{1 + \exp\!\left[500\,(x_1 + x_2)\right]}, & \mathbf{x} \in \Omega, \\[8pt]
u(\mathbf{x}, t) = \dfrac{1}{1 + \exp\!\left[500\,(x_1 + x_2 - t)\right]}, & (\mathbf{x}, t) \in \partial\Omega \times (0, T],
\end{cases}\label{eq:burgers}
\end{equation}
where $\Omega = (-1, 1)^2$, $T = 1$. The exact solution is
\begin{equation}
u(\mathbf{x}, t) = \frac{1}{1 + \exp\!\left[500\,(x_1 + x_2 - t)\right]}.\label{eq:burgers_exact}
\end{equation}

The solution develops a sharp internal layer along the diagonal $x_1 + x_2 = t$, which sweeps across the domain as $t$ increases from $0$ to $1$. The training configuration and validation grid follow exactly those of the previous problem, with 200,000 validation points evaluated on a $200 \times 200$ spatial grid over $(-1,1)^2$ at 5 uniformly spaced time levels.

\begin{table}[ht!]
\centering
\caption{{Nonlinear problem \eqref{eq:burgers}}: minimum $\mathcal{E}_2$, $\mathcal{E}_\infty$, and ${\mathcal{E}_{\mathrm{mse}}}$ at each adaptive stage.}\label{tab:burgers}
\resizebox{\textwidth}{!}{%
\begin{tabular}{lccccc}
\toprule
& \textbf{initial stage} & \textbf{adaptive stage 1} & \textbf{adaptive stage 2} & \textbf{adaptive stage 3} & \textbf{adaptive stage 4} \\
\midrule
Iter range & $1$--$8000$ & $8001$--$24000$ & $24001$--$48000$ & $48001$--$80000$ & $80001$--$120000$ \\
$\mathcal{E}_2$ min & $1.20\times10^{-1}$ & $1.83\times10^{-2}$ & $5.83\times10^{-3}$ & $3.59\times10^{-3}$ & $1.70\times10^{-3}$ \\
$\mathcal{E}_\infty$ min & $4.85\times10^{-1}$ & $3.10\times10^{-1}$ & $1.91\times10^{-1}$ & $1.00\times10^{-1}$ & $3.81\times10^{-2}$ \\
${\mathcal{E}_{\mathrm{mse}}} $ min  & $1.00\times10^{-2}$ & $2.35\times10^{-4}$ & $2.38\times10^{-5}$ & $9.04\times10^{-6}$ & $2.03\times10^{-6}$ \\
\bottomrule
\end{tabular}}
\end{table}

\begin{figure}[tbhp!]
\centering
\begin{subfigure}{\textwidth}
\centering
\includegraphics[width=0.87\textwidth]{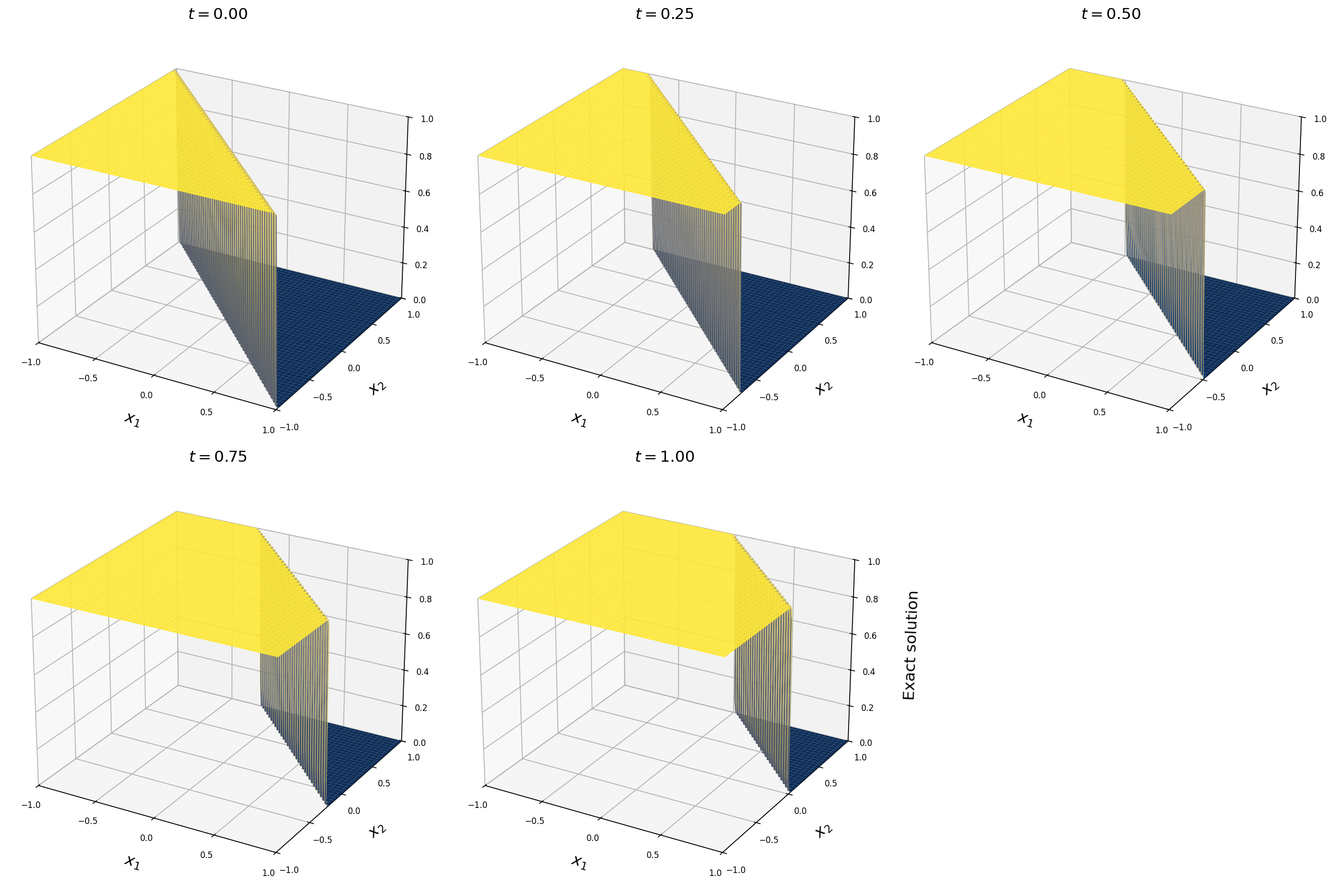}
\caption{Exact solution}
\end{subfigure}
\vspace{0.5cm}
\begin{subfigure}{\textwidth}
\centering
\includegraphics[width=0.87\textwidth]{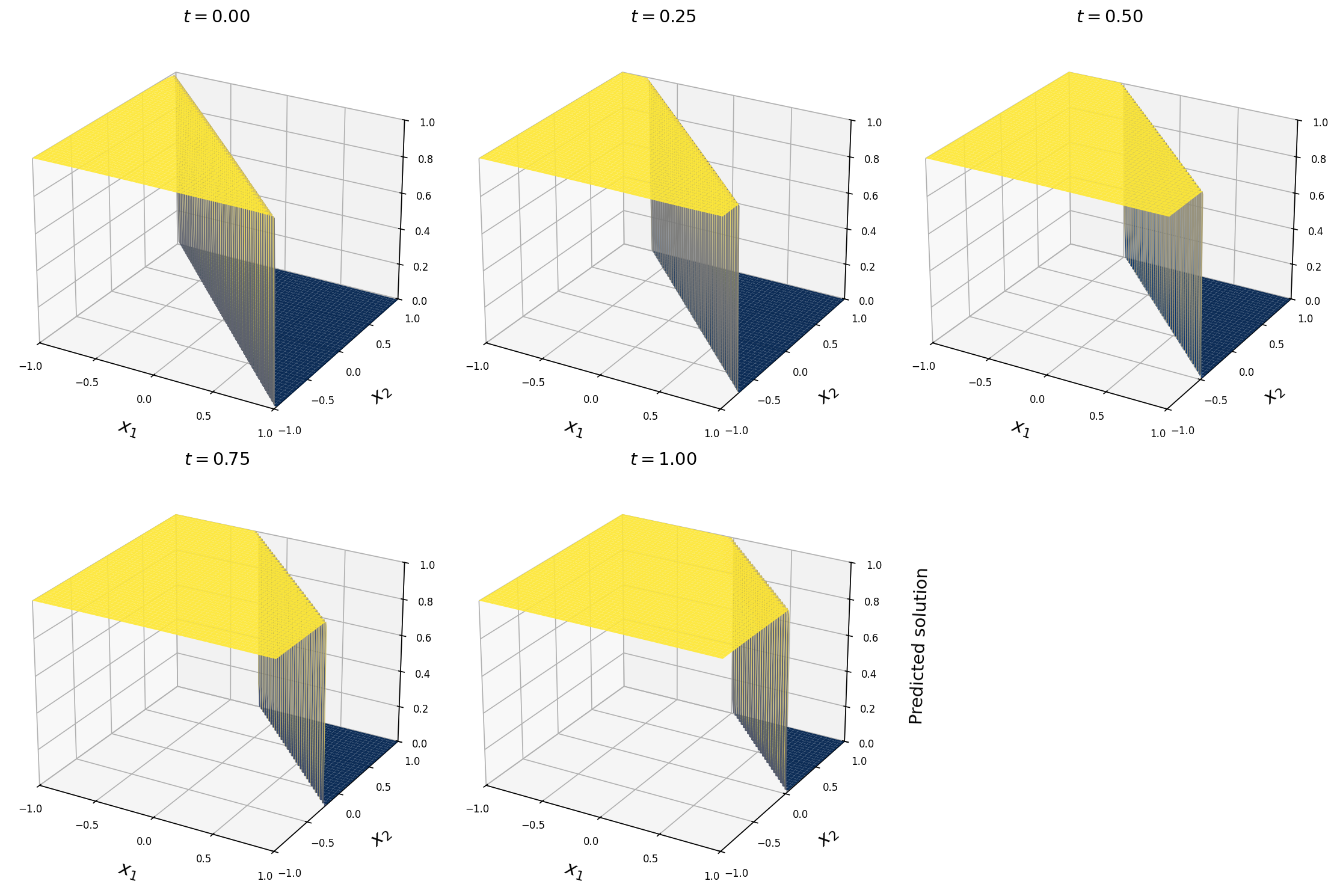}
\caption{DAS-PINNs predicted solution}
\end{subfigure}
\caption{{Nonlinear problem  \eqref{eq:burgers}}:  exact solution (top) and predicted solution (bottom) at five time levels $t \in \{0, 0.25, 0.50, 0.75, 1.00\}$.}\label{fig:prob3_solution}
\end{figure}

\begin{figure}[ht!]
\centering
\includegraphics[width=0.9\linewidth]{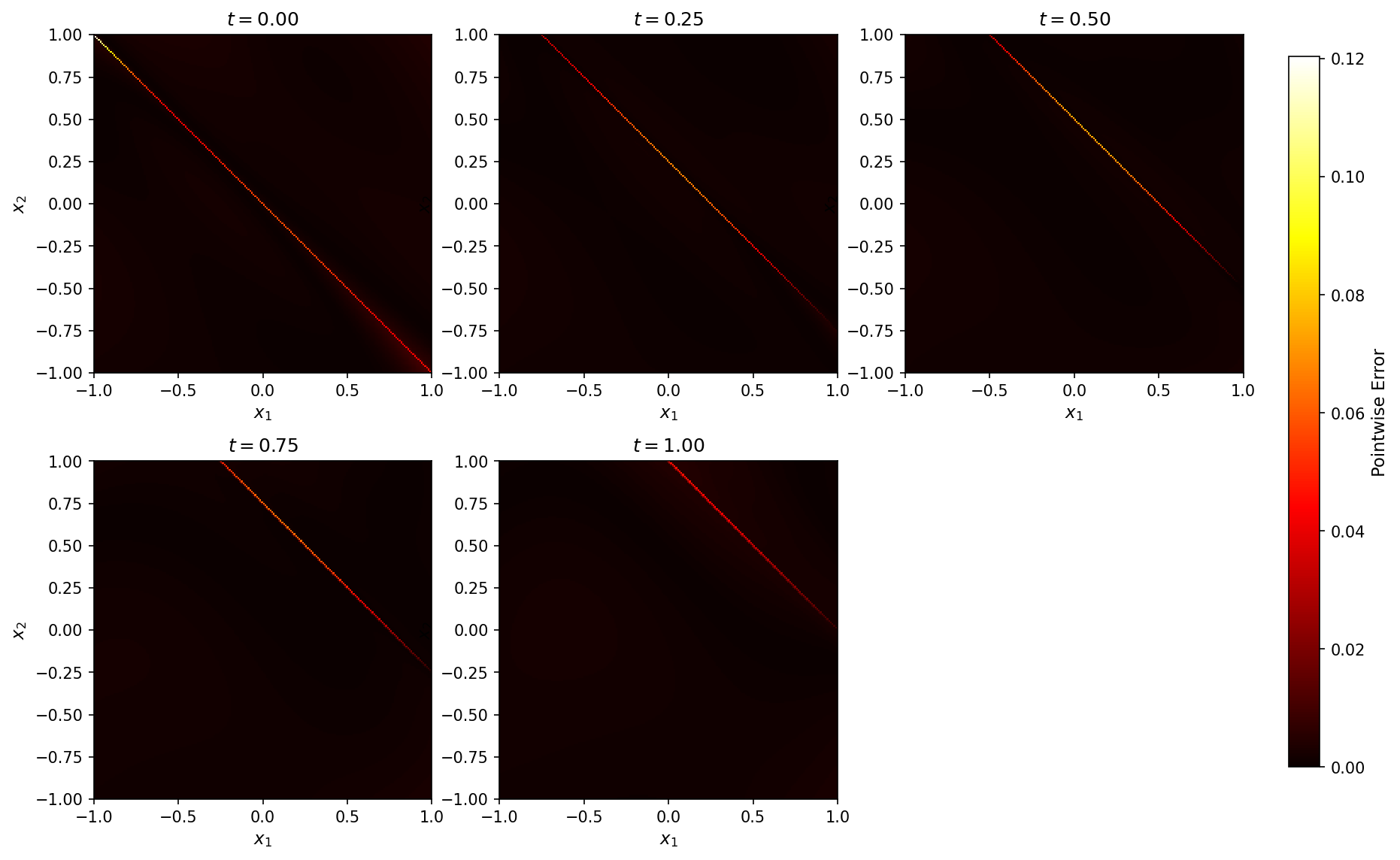}
\caption{{Nonlinear problem \eqref{eq:burgers}}:  pointwise error at five time levels $t \in \{0, 0.25, 0.50, 0.75, 1.00\}$.}\label{fig:prob3_pointwise_error}
\end{figure}
	
\begin{figure}[ht!]
		\centering
		\includegraphics[width=0.85\linewidth]{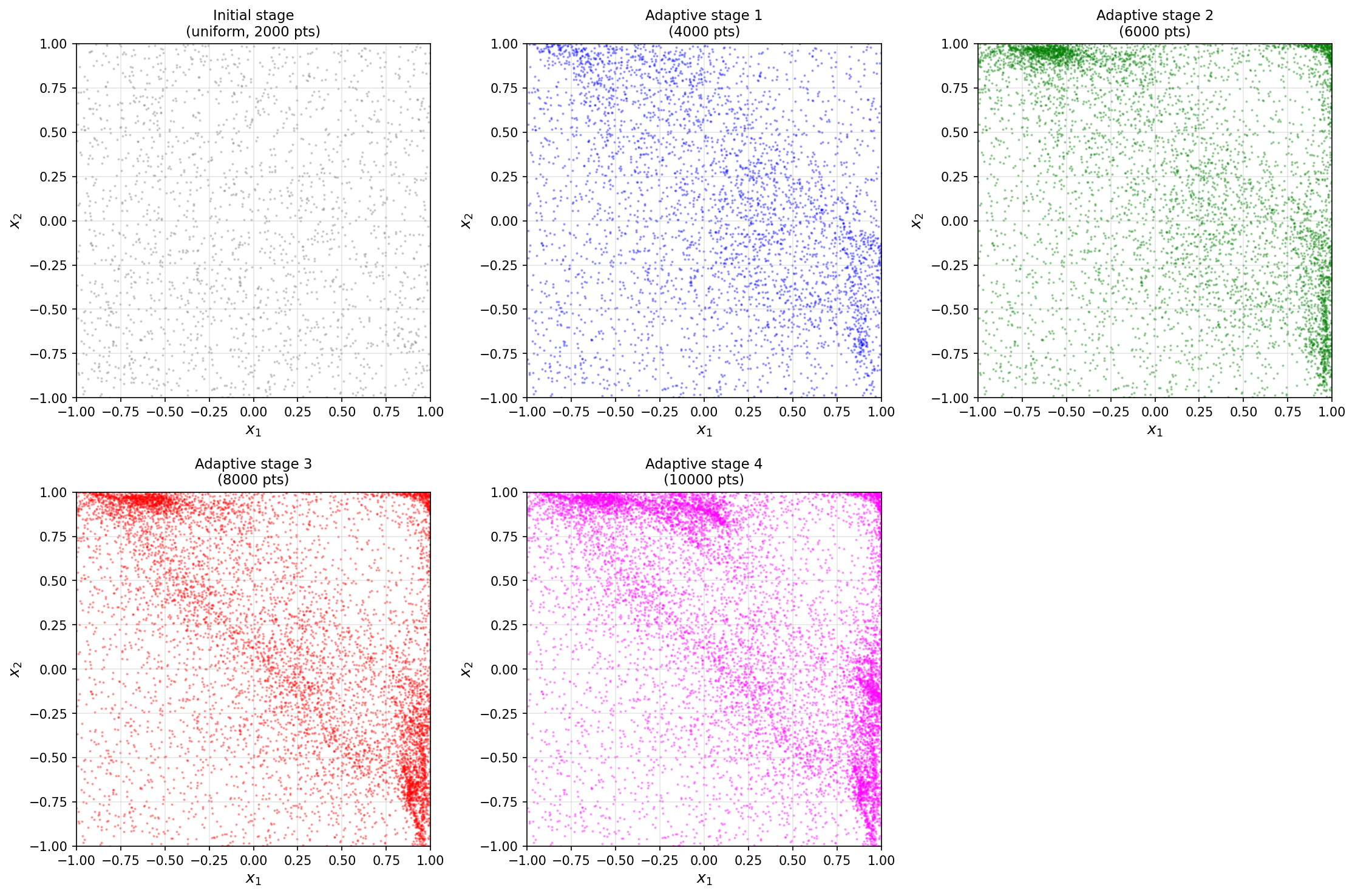}
	\caption{{Nonlinear problem \eqref{eq:burgers}}:  spatial distribution of collocation 
		points at each stage of DAS-PINNs.}\label{fig:prob3_spatial_points}
	\end{figure}

Table~\ref{tab:burgers} summarises the minimum $\mathcal{E}_2$, $\mathcal{E}_\infty$, and ${\mathcal{E}_{\mathrm{mse}}}$ achieved at each stage. All three metrics decrease consistently as the adaptive stages progress, with the most significant reduction occurring after the  first adaptive stage. By the final stage, $\mathcal{E}_2$, $\mathcal{E}_\infty$, and ${\mathcal{E}_{\mathrm{mse}}}$ reach $1.70\times10^{-3}$, $3.81\times10^{-2}$, and $2.03\times10^{-6}$ respectively, demonstrating that DAS-PINNs successfully resolves the sharp internal layer despite its extreme steepness. Figure~\ref{fig:prob3_solution} shows the exact and DAS-{PINNs} predicted solutions at {five} time levels. The predicted solution faithfully reproduces the sharp diagonal front at each time level, with the transition from one to zero captured cleanly across the entire time interval. Figure~\ref{fig:prob3_pointwise_error} shows the pointwise error at the same time levels. The error is almost entirely confined to a thin band along the front $x_1 + x_2 = t$, while remaining negligibly small throughout the rest of the domain, which reflects the inherent difficulty of resolving such an extremely steep gradient with a finite number of collocation points. Figure~\ref{fig:prob3_spatial_points} shows the spatial distribution of collocation points across all five stages. Starting from a uniform distribution in the initial stage, DAS-PINNs progressively concentrates points toward the upper-right corner of the domain, where the sharp layer intersects the boundary. As training advances, the point distribution spreads along the full diagonal trajectory of the front, indicating that KRnet is tracking the high-residual region as it sweeps across the domain.
	
\subsection{High-dimensional problems}\label{sec:hd}
In this section, we consider two time-dependent problems in high-dimensional spatial domains. The first is a parabolic problem with a localised Gaussian source, solved for $d = 6$ and $d = 8$. The second is a high-dimensional advection equation, solved for $d = 6$. In all cases, the full spatiotemporal domain $\mathbb{R}^d \times [0, T]$ is treated within the DAS-PINNs framework, and the performance of DAS-PINNs is compared against uniform sampling to demonstrate the benefit of adaptive collocation in high dimensions.

\subsubsection{A high-dimensional parabolic problem}\label{sec:hdpp}
We consider the following parabolic {PDE} problem:
\begin{equation}
	\begin{cases}
		\dfrac{\partial u}{\partial t} - \Delta u = f(\mathbf{x}, t), 
		& (\mathbf{x}, t) \in \Omega \times (0, T], \\[8pt]
		u(\mathbf{x}, 0) = 0, & \mathbf{x} \in \Omega, \\[8pt]
		u(\mathbf{x}, t) = {0}, 
		& (\mathbf{x}, t) \in \partial\Omega \times (0, T],
	\end{cases}
	\label{eq:parabolic_hd}
\end{equation}
where $\Omega = [-1, 1]^d$ and $T = 1$. The exact solution is {given by}
\begin{equation}
	u(\mathbf{x}, t) = \exp\!\left(-10\|\mathbf{x}\|^2\right)\left(1 - e^{-5t}\right),
	\label{eq:parabolic_hd_exact}
\end{equation}
and the source term $f(\mathbf{x}, t)$ is derived analytically as
\begin{equation}
	f(\mathbf{x}, t) = \exp\!\left(-10\|\mathbf{x}\|^2\right)
	\left[5e^{-5t} + \left(20d - 400\|\mathbf{x}\|^2\right)
	\left(1 - e^{-5t}\right)\right],
	\label{eq:parabolic_hd_source}
\end{equation}
where $\|\mathbf{x}\|^2 = \sum_{i=1}^{d} x_i^2$. The solution evolves from zero at $t = 0$ toward a steady-state Gaussian profile centred at the origin, with the rate of growth controlled by the factor $\left(1 - e^{-5t}\right)$. The Gaussian profile becomes 
increasingly concentrated near the origin as $d$ increases, since the volume of the region where $\exp(-10\|\mathbf{x}\|^2)$ is appreciable shrinks rapidly with dimension, making this a genuinely challenging test for uniform sampling strategies.
	
The solution network uses 6 fully connected layers with 64 hidden neurons per layer and a $\tanh$ activation function. For KRnet, 6 affine coupling layers are used, each consisting of two fully connected layers with 64 neurons per layer. Both networks are 
trained with $n_{\text{train}} = 10{,}000$ collocation points and batch size 5000, over one initial stage followed by four adaptive stages, giving a total of 90,000 training iterations, accumulating 50,000 collocation points by the final stage. The uniform 
sampling baseline uses 50,000 collocation points trained in a single stage with the same total iteration budget, making the comparison fair in terms of both total points and total iterations. The model is evaluated on a near-origin validation grid over $[-0.1, 0.1]^d$ at 5 uniformly spaced time levels, comprising $5^6 \times 5 = 78{,}125$ points for $d = 6$ and $3^8 \times 5 = 32{,}805$ points for $d = 8$.
In preliminary experiments, the solution network struggled to satisfy the zero initial condition $u(\mathbf{x}, 0) = 0$ accurately, particularly in early training stages where the residual loss dominated. To address this, we impose a hard constraint on the network output by defining $u_\theta = t \cdot y_\theta$, where $y_\theta$ is the raw network output. This transformation ensures that the initial condition $u(\mathbf{x}, 0) = 0$ is satisfied exactly by construction, removing it from the loss function and allowing the network to focus entirely on fitting the PDE residual and boundary conditions.

\begin{figure}[ht!]
		\centering
		\begin{subfigure}{\textwidth}
			\centering
			\includegraphics[width=\textwidth]{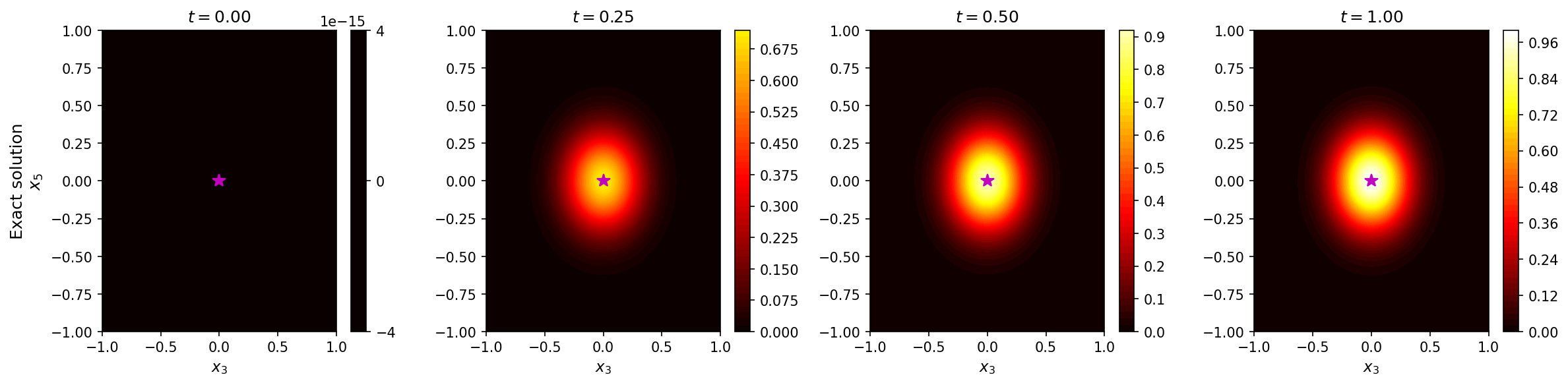}
		\end{subfigure}
		\vspace{0.5cm}
		\begin{subfigure}{\textwidth}
			\centering
			\includegraphics[width=\textwidth]{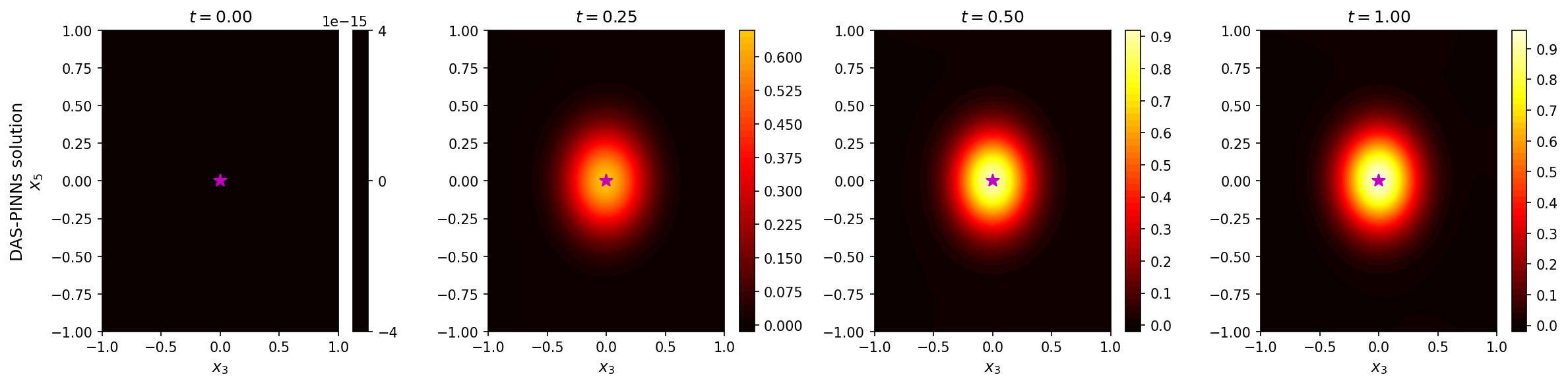}
		\end{subfigure}
		\vspace{0.5cm}
		\begin{subfigure}{\textwidth}
			\centering
			\includegraphics[width=\textwidth]{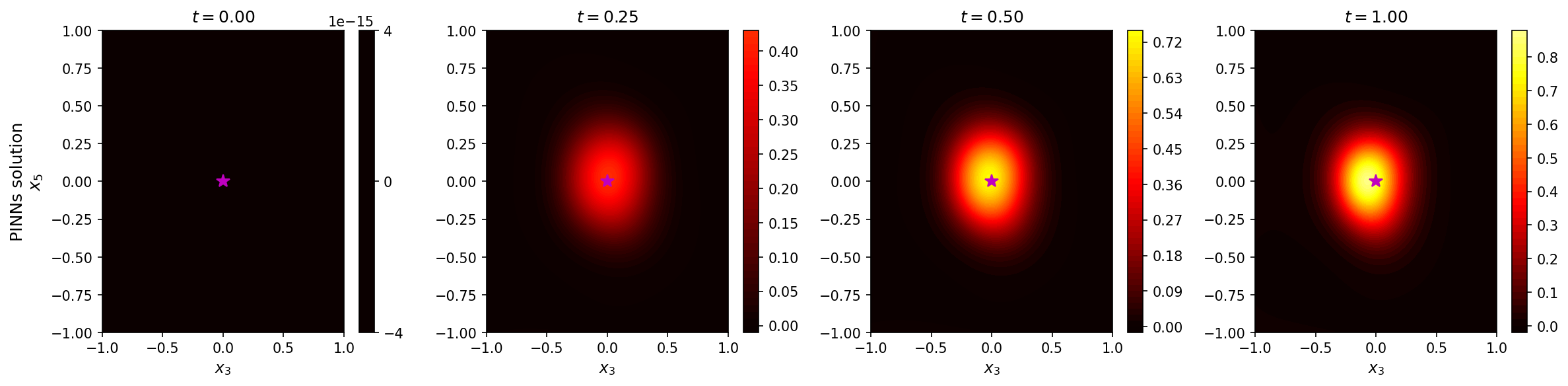}
		\end{subfigure}
		\caption{{High-dimensional parabolic problem \eqref{eq:parabolic_hd}} $d = 6$: exact solution (top row), 
			DAS-PINNs predicted solution (middle row), and PINNs with uniform sampling 
			(bottom row) at {four time levels $t \in \{0, 0.25, 0.50, 1.00\}$}, 
			shown as a 2D slice along $(x_3, x_5)$.}
		\label{fig:prob4_dim7_solution}
\end{figure}

\begin{table}[ht!]
	\caption{{High-dimensional parabolic problem \eqref{eq:parabolic_hd}}: comparison of 
	DAS-PINNs and uniform sampling.}
	\centering
	\begin{tabular}{l cccc}
		\toprule
		& \multicolumn{2}{c}{$d = 7$ (spatial dim = 6)} & \multicolumn{2}{c}{$d = 9$ (spatial dim = 8)} \\
		\cmidrule(lr){2-3} \cmidrule(lr){4-5}
		& DAS-PINNs & Uniform & DAS-PINNs & Uniform \\
		\midrule
		$	\mathcal{E}_{\mathrm{mse}}$ & $\mathbf{3.83 \times 10^{-4}}$ & $2.12 \times 10^{-2}$ & $\mathbf{1.81 \times 10^{-3}}$ & $2.23 \times 10^{-1}$ \\
		$\mathcal{E}_2$  & $\mathbf{3.22 \times 10^{-2}}$ & $2.40 \times 10^{-1}$ & $\mathbf{8.77 \times 10^{-2}}$ & $9.74 \times 10^{-1}$ \\
		\bottomrule
	\end{tabular}\label{tab:parabolic_hd}
\end{table}

Table~\ref{tab:parabolic_hd} summarises the {errors} ${\mathcal{E}_{\mathrm{mse}}}$ and $\mathcal{E}_2$ achieved by DAS-PINNs and uniform sampling for spatial dimensions $d = 6$ and $d = 8$. In both cases, DAS-PINNs outperforms uniform sampling by a substantial margin, with the gap widening as the spatial dimension increases. At $d = 8$, uniform sampling fails to resolve the localised Gaussian structure with $\mathcal{E}_2$ approaching unity, while DAS-PINNs maintains reasonable accuracy. This progressive deterioration of uniform sampling with increasing spatial dimension is a direct manifestation of the curse of dimensionality, and demonstrates that adaptive collocation is not merely beneficial but essential for resolving localised  solutions in high-dimensional spatiotemporal domains.

\begin{figure}[tbhp!]
	\centering
	\begin{subfigure}{\textwidth}
		\centering
		\includegraphics[width=\textwidth]{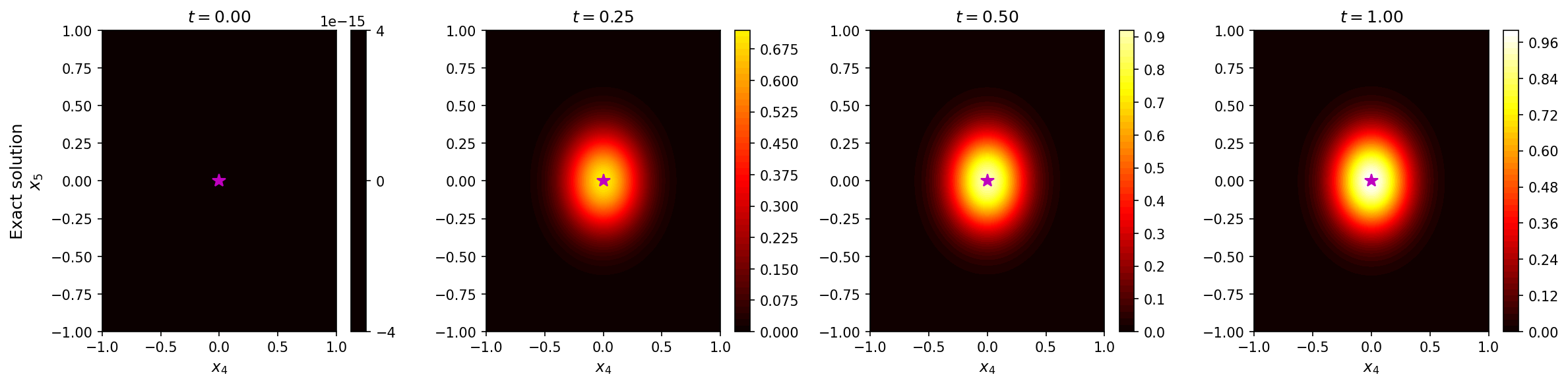}
	\end{subfigure}
	\vspace{0.5cm}
	\begin{subfigure}{\textwidth}
		\centering
		\includegraphics[width=\textwidth]{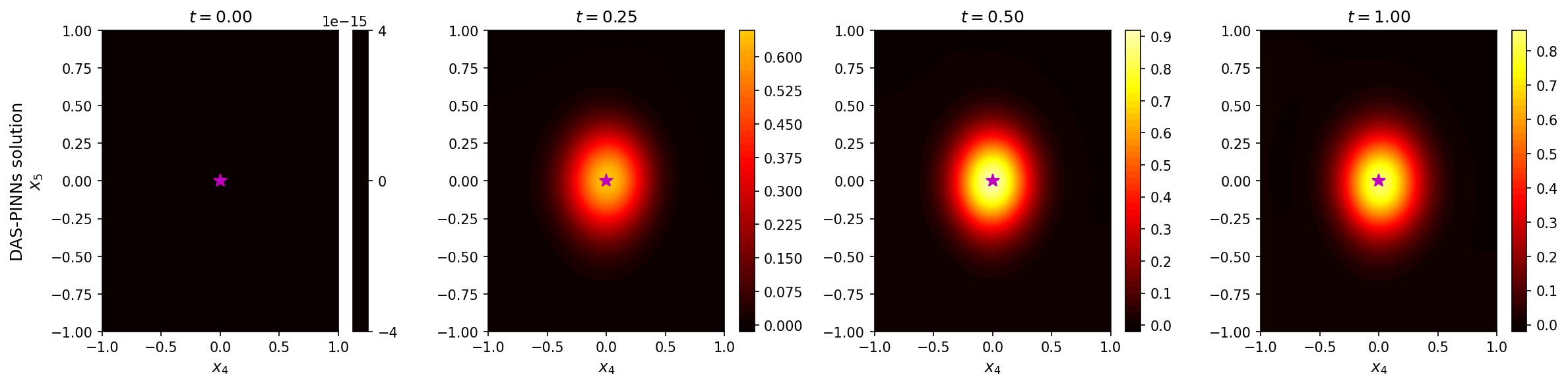}
	\end{subfigure}
	\vspace{0.5cm}
	\begin{subfigure}{\textwidth}
		\centering
		\includegraphics[width=\textwidth]{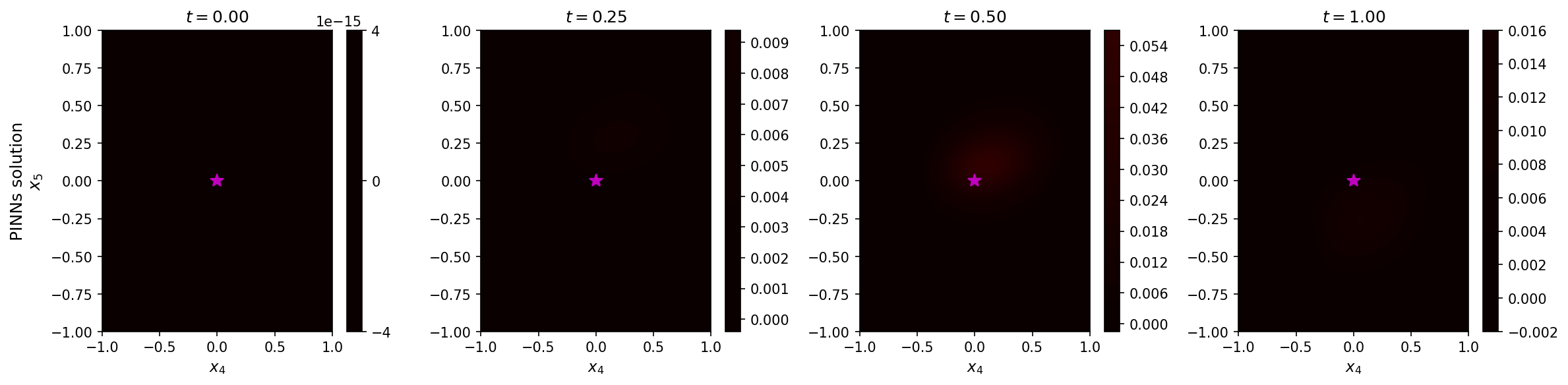}
	\end{subfigure}
	\caption{{High-dimensional parabolic problem \eqref{eq:parabolic_hd}} $d = 8$: exact solution (top row), 
		DAS-PINNs predicted solution (middle row), and PINNs with uniform sampling 
		(bottom row) at {four time levels $t \in \{0, 0.25, 0.50, 1.00\}$},  
		shown as a 2D slice along $(x_4, x_5)$.}
	\label{fig:prob4_dim9_solution_x4_x5}
\end{figure}

 Figure~\ref{fig:prob4_dim7_solution} shows the exact solution alongside the predictions from DAS-PINNs and uniform sampling at five time levels, visualised as a 2D slice along $(x_3, x_5)$, with the remaining dimensions fixed at zero to pass through the peak of the Gaussian profile. At $t = 0$ the solution is identically zero since the exact solution contains the factor $(1 - e^{-5t})$, which vanishes at $t = 0$; the star marker serves only as a reference for the peak location from which the Gaussian profile emerges for $t > 0$. The DAS-PINNs prediction closely matches the exact solution across all time levels, capturing both the shape and peak amplitude as the solution grows from zero toward its steady state. Uniform sampling noticeably underestimates the peak throughout---the compressed colorbars in the bottom row make this immediately apparent, which is a direct consequence of too few points landing near the origin in six dimensions. Slices along other dimension pairs were also examined and showed no significantly different behaviour, confirming that the results are representative of the network's accuracy across the full spatial domain. Figure~\ref{fig:prob4_dim9_solution_x4_x5} shows the corresponding results for $d = 8$, visualised as a 2D slice along $(x_4, x_5)$. DAS-PINNs continues to capture the Gaussian profile accurately throughout the time interval. Uniform sampling, by contrast, produces a response that is nearly zero at all time levels, the colorbars are two orders of magnitude smaller than the exact solution, indicating a complete failure to resolve the localised structure in the nine-dimensional spatiotemporal domain. Moving from spatial dimension, $d = 6$ to $d = 8$, uniform sampling goes from poor to essentially	useless, a clear sign that the curse of dimensionality is at work. DAS-PINNs, by focusing points where the solution actually lives, avoids this trap entirely. 
 
 \begin{figure}[tbhp!]
	\centering
	\includegraphics[width=\linewidth]{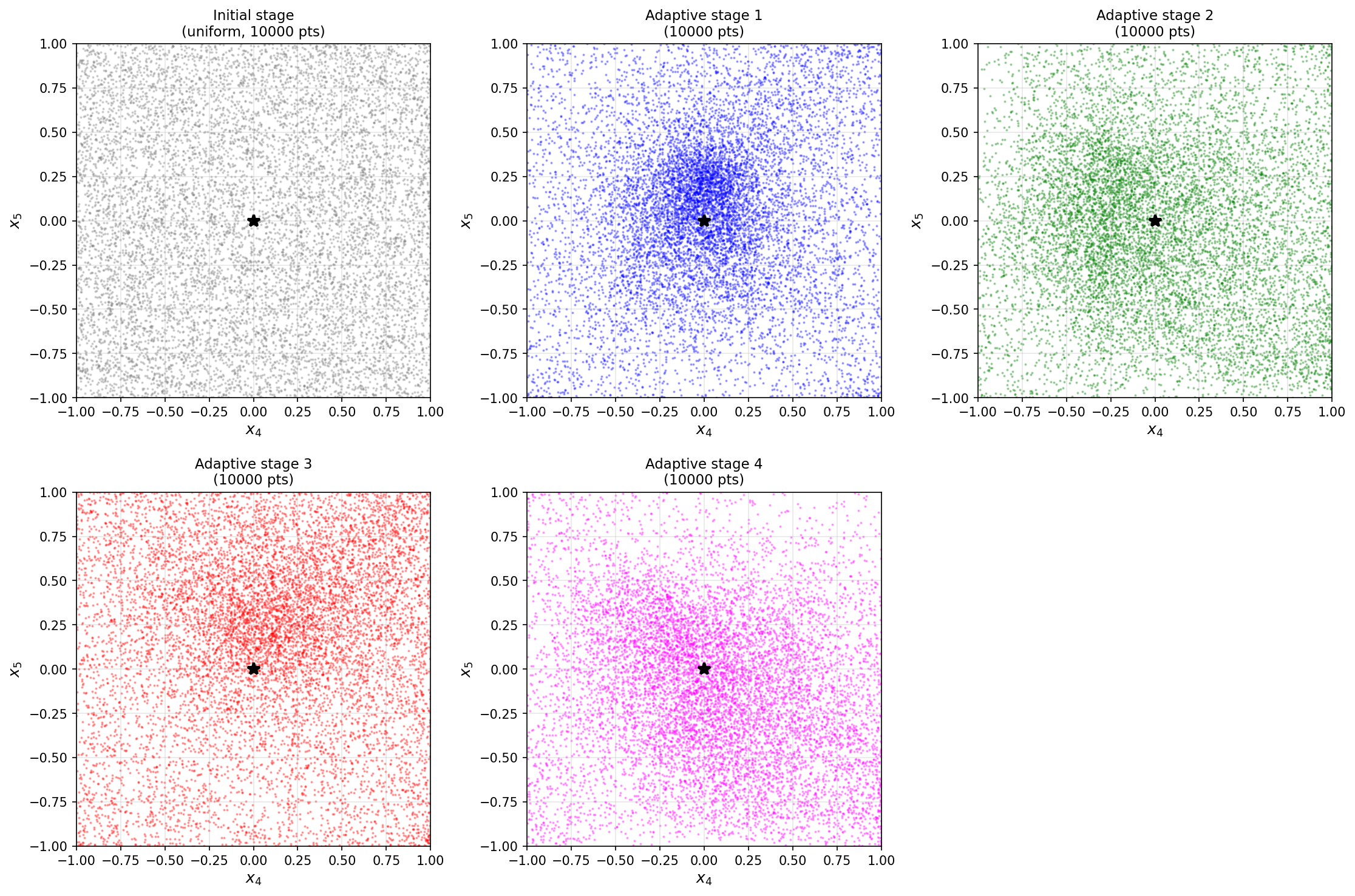}
	\caption{{High-dimensional parabolic problem \eqref{eq:parabolic_hd}} $d = 8$: spatial distribution of 
		newly added collocation points at each stage of DAS-PINNs, projected onto 
		$(x_4, x_5)$. The star marks the peak location at the origin.}\label{fig:prob4_dim9_spatial_points}
\end{figure}
 
 Figure~\ref{fig:prob4_dim9_spatial_points} shows the newly added collocation points at each stage, projected onto the $(x_4, x_5)$ plane; each panel shows only the points added at that stage, not the full accumulated training set. Starting from a uniform distribution, KRnet immediately identifies the high-residual region near the origin and concentrates the first batch of adaptive points there. As training progresses and the approximation near the origin improves, the residual spreads more evenly and KRnet responds by distributing subsequent points more  broadly. The sampler is not targeting a fixed region but reacting to where the network is currently weakest, which shifts outward as the localised structure near the origin is resolved. Projections onto other dimension pairs were also examined and showed no significantly different behaviour.

\subsubsection{{A high-dimensional hyperbolic problem}}\label{High-Dimensional Advection Equation}
The second high-dimensional test problem is a pure transport equation in which a narrow Gaussian profile propagates along the main diagonal of the domain{:}	
\begin{equation}
		\begin{cases}
			\dfrac{\partial u}{\partial t} = -\displaystyle\sum_{i=1}^{d} 
			\dfrac{\partial u}{\partial x_i}, 
			& (\mathbf{x}, t) \in \Omega \times (0, T], \\[8pt]
			u(\mathbf{x}, 0) = \exp\!\left(-10\sum_{i=1}^{d} x_i^2\right), 
			& \mathbf{x} \in \Omega, \\[8pt]
			u(\mathbf{x}, t) = \exp\!\left(-10\sum_{i=1}^{d} (x_i - t)^2\right), 
			& (\mathbf{x}, t) \in \partial\Omega \times (0, T],
		\end{cases}
		\label{eq.hdhp}
	\end{equation}
where $\Omega = (-0.2, 1.2)^6$ and $T = 1$. The exact solution is {given by}
	\begin{equation}
		u(\mathbf{x}, t) = \exp\!\left(-10\sum_{i=1}^{d}(x_i - t)^2\right).
		\label{eq:advection_hd_exact}
	\end{equation}
The peak centre travels along the trajectory $(t, t, \ldots, t)$, moving from the origin at $t = 0$ to the corner $(1, 1, \ldots, 1)$ at $t = 1$. With a sharpness coefficient of 10 in the exponent, the Gaussian profile remains highly localised throughout its motion. Unlike the steady-state setting, where the high-residual region is fixed, the adaptive sampler here must track a moving target in a six-dimensional spatial domain---the peak location shifts continuously with time, requiring KRnet to identify and follow the high-residual region across the full spatiotemporal domain $(-0.2, 1.2)^6 \times [0, 1]$.

To handle this spatiotemporal localisation, the training configuration follows that of the previous problem, except that $n_{\text{train}} = 5{,}000$ collocation points are used with batch size 1000. The model is evaluated on a tensor product validation grid of $5^6 \times 5 = 78{,}125$ points over $(-0.2, 1.2)^6$ at 5 uniformly spaced time levels.

\begin{table}[ht!]
	\centering
	\caption{{High-dimensional hyperbolic problem \eqref{eq.hdhp}}: minimum $\mathcal{E}_2$, 
		$\mathcal{E}_\infty$, and ${\mathcal{E}_{\mathrm{mse}}}$ at each adaptive stage.}
	\label{tab:advection_hd}
	\resizebox{\textwidth}{!}{%
		\begin{tabular}{lccccc}
			\toprule
			& \textbf{initial stage} & \textbf{Adaptive stage 1} 
			& \textbf{Adaptive stage 2} & \textbf{Adaptive stage 3} & \textbf{Adaptive stage 4} \\
			\midrule
			Iter range 
			& $1$--$15000$ & $15001$--$45000$ 
			& $45001$--$90000$ & $90001$--$150000$ & $150001$--$225000$ \\
			$\mathcal{E}_2$ min 
			& $9.51\times10^{-1}$ & $1.31\times10^{-1}$ 
			& $8.24\times10^{-2}$ & $6.35\times10^{-2}$ & $4.57\times10^{-2}$ \\
			$\mathcal{E}_\infty$ min 
			& $9.83\times10^{-1}$ & $2.28\times10^{-1}$ 
			& $1.59\times10^{-1}$ & $1.29\times10^{-1}$ & $9.12\times10^{-2}$ \\
			${\mathcal{E}_{\mathrm{mse}}}$ min  
			& $5.06\times10^{-1}$ & $9.56\times10^{-3}$ 
			& $3.81\times10^{-3}$ & $2.26\times10^{-3}$ & $1.17\times10^{-3}$ \\
			\bottomrule
	\end{tabular}}
\end{table}

Table~\ref{tab:advection_hd} summarises the minimum $\mathcal{E}_2$, $\mathcal{E}_\infty$, and ${\mathcal{E}_{\mathrm{mse}}}$ at each stage. All three metrics decrease steadily across the adaptive stages, with the most substantial drop occurring after the first adaptive stage. The continued improvement through subsequent stages demonstrates that KRnet is successfully tracking the moving Gaussian peak across the spatiotemporal domain, progressively concentrating collocation points along its trajectory as training advances. 

\begin{figure}[tbhp!]
	\centering
	\includegraphics[width=0.95\linewidth]{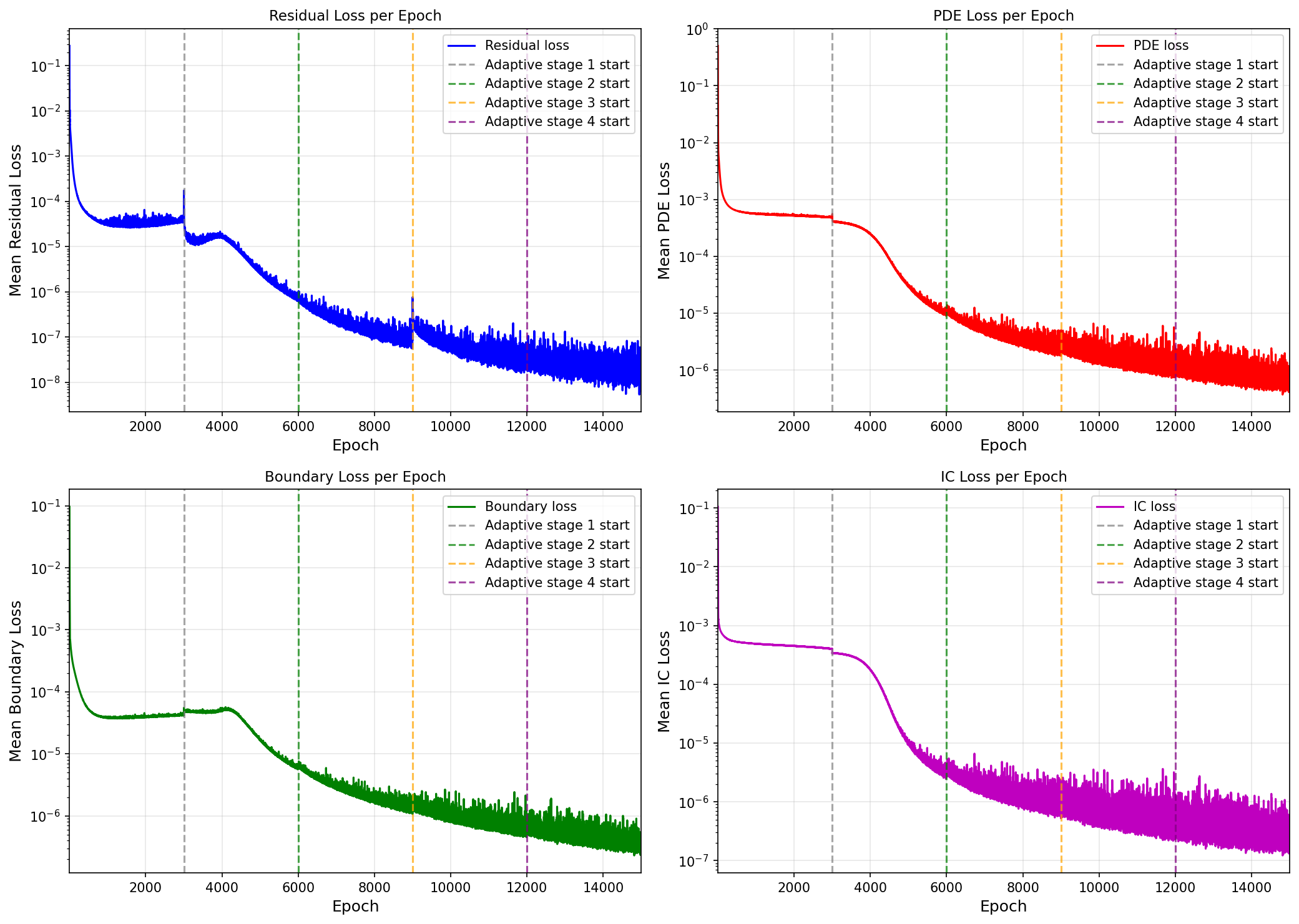}
	\caption{{High-dimensional hyperbolic problem \eqref{eq.hdhp}}): evolution of the loss components---residual loss, PDE loss, boundary loss, and initial condition loss---per epoch. Dashed vertical lines mark the start of each adaptive stage.}\label{fig:prob5_loss}
\end{figure}

\begin{figure}[tbhp!]
\centering
\includegraphics[width=0.95\linewidth]{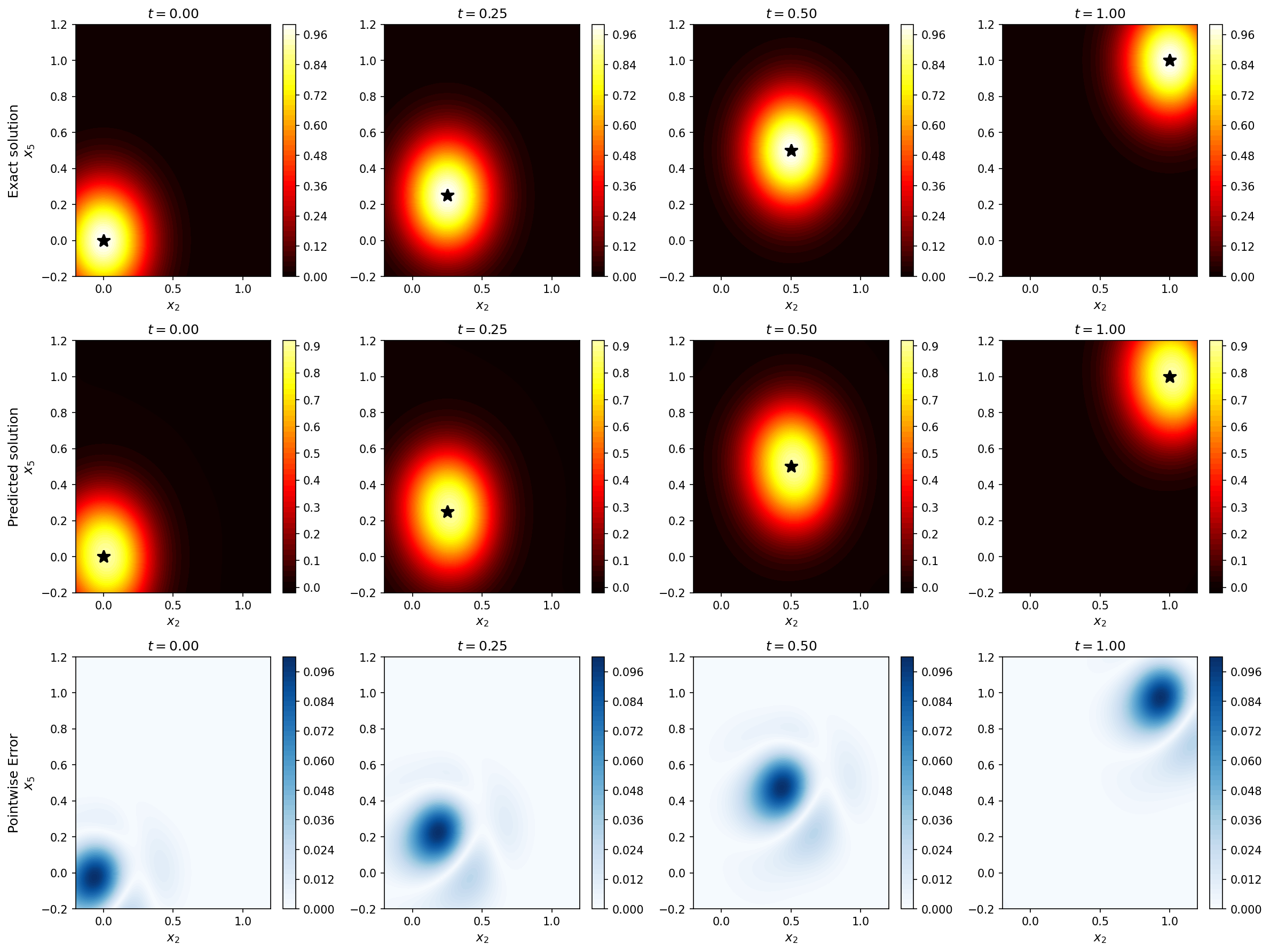}
\caption{{High-dimensional hyperbolic problem \eqref{eq.hdhp}}): exact solution (top row), predicted solution (middle row), and pointwise error (bottom row) at {four time levels $t \in \{0, 0.25, 0.50, 1.00\}$}, 
shown as a 2D slice along $(x_2, x_5)$.}\label{fig:prob5_solution_slices}
\end{figure}

Figure~\ref{fig:prob5_loss} shows the evolution of all four loss components across training. All components decrease sharply during the initial stage, with a brief transient increase at the start of adaptive stage 1 as the collocation points are resampled. The most significant reduction occurs between adaptive stages 1 and 2, after which all four components continue to decay smoothly and consistently through the remaining stages, reaching values of order $10^{-6}$ to $10^{-7}$ by the end of training. The boundary and initial condition losses decrease faster than the residual loss, reflecting the relative difficulty of satisfying the PDE in the interior of a seven-dimensional spatiotemporal domain $(-0.2, 1.2)^6 \times [0, 1]$.

Figure~\ref{fig:prob5_solution_slices} shows the exact solution, predicted solution, and pointwise error at five time levels as a 2D slice along $(x_2, x_5)$, with all remaining dimensions fixed at the peak location $x_k = t$. The star marker  tracks the peak centre as it travels along the diagonal of the domain. The predicted solution closely matches the exact solution across all time levels, capturing both the shape and magnitude of the Gaussian profile throughout its motion. The pointwise error remains localised near the peak centre and stays small relative to the solution amplitude, demonstrating that DAS-PINNs successfully resolves the moving localised structure in the full seven-dimensional spatiotemporal domain. Since the exact solution is symmetric across all spatial dimensions, this slice is representative of the network's accuracy throughout the entire domain.

\section{Conclusion}\label{sec:Conclusion}
	
In this work, the DAS-PINNs  methodology is applied  to time-dependent high-dimensional PDEs, demonstrating that treating the spatiotemporal domain as a unified high-dimensional space provides a natural and effective setting for deep adaptive sampling. By training KRnet to approximate the residual distribution across both space and time, collocation  points are automatically concentrated where the solution was most difficult to learn, without any explicit time-marching or problem-specific intervention. Numerical experiments across five benchmark problems confirmed this behaviour, with KRnet consistently identifying and adapting to high-residual regions as training progressed. In low-dimensional settings, DAS-PINNs successfully resolved a range of challenging time-dependent problems with moving and localised solution features. In high-dimensional settings, the advantage of adaptive sampling becomes increasingly apparent as the spatial dimension grows. 
These results confirm that deep adaptive sampling is an effective approach for time-dependent high-dimensional PDEs, opening the way for further exploration in more complex settings such as fractional and nonlocal partial differential equations. 
	
\section*{Acknowledgments}
This work was supported by EPSRC grant EP/Y028783/1. 

\bibliographystyle{elsarticle-num}
\bibliography{mref}
\end{document}